\input amstex
\magnification=\magstep1 
\baselineskip=13pt
\documentstyle{amsppt}
\vsize=8.7truein \CenteredTagsOnSplits \NoRunningHeads
\def\today{\ifcase\month\or
  January\or February\or March\or April\or May\or June\or
  July\or August\or September\or October\or November\or December\fi
  \space\number\day, \number\year}
 
 \def\EE{\bold{E}\thinspace}
 \def\PP{\bold{P}\thinspace}

 \def\spa{\operatorname{span}}
 \def\sgn{\operatorname{sign}}
 \def\cov{\bold{cov}\thinspace}
  \def\var{\bold{var}\thinspace}
  \def\dist{\operatorname{dist}}
 
 \topmatter

\title  The number of  graphs and a random graph with a given degree sequence \endtitle
\author Alexander Barvinok and J.A. Hartigan 
\endauthor
\address Department of Mathematics, University of Michigan, Ann Arbor,
MI 48109-1043, USA
\endaddress
\email barvinok$\@$umich.edu \endemail
\address Department of Statistics, Yale University, New Haven, CT 06520-8290 \endaddress
\email john.hartigan$\@$yale.edu \endemail 
\thanks The research of the first author was partially supported by NSF Grant DMS 0856640 and a
United States - Israel BSF grant 2006377. \endthanks
\date November 2011
\enddate
\keywords graphs, degree sequences, asymptotic formulas
\endkeywords
\subjclass 05A16, 05C07, 05C30, 52B55, 60F05
 \endsubjclass
\abstract We consider the set of all graphs on $n$ labeled vertices with prescribed degrees
$D=\left(d_1, \ldots, d_n\right)$. For a wide class of {\it tame} degree sequences $D$ we 
obtain a computationally efficient asymptotic  formula approximating the number of graphs within 
a relative error which approaches 0 as $n$ grows. As a corollary, we prove that the structure of  a random graph with a given tame degree sequence $D$ is well described by a certain 
{\it maximum entropy matrix} computed from $D$.
We also establish an asymptotic formula for the 
number of bipartite graphs with prescribed degrees of vertices, or, equivalently, for the number of 0-1 matrices
with prescribed row and column sums.
\endabstract
\endtopmatter

\document

\head 1. Introduction and main results \endhead

\subhead (1.1) Graphs and their degree sequences \endsubhead
Let $D=\left(d_1, \ldots, d_n\right)$ be a vector of positive integers and let $G(D)$ be the 
set of all graphs (undirected, with no loops or multiple edges) on the set 
$\{1, \ldots, n\}$ of vertices such that the degree of the $k$-th vertex is $d_k$ for 
$k=1, \ldots, n$. Equivalently, $G(D)$ is the set of all $n \times n$ symmetric matrices 
with 0-1 entries, zero trace and row (column) sums $d_1, \ldots, d_n$. 
We assume that 
$$d_1 + \ldots + d_n \equiv 0 \mod 2, \tag1.1.1$$
since otherwise the set $G(D)$ is empty.

The theorem of Erd\H{o}s and Gallai, see, for example, Theorem 6.3.6 of  \cite{BR91}, states the necessary and sufficient conditions for the 
existence of a graph with the given degree sequence. Without loss of generality, we assume that 
$$d_1 \geq d_2 \geq \ldots \ \geq d_n.$$ 
Then, the necessary and sufficient condition for $G(D)$ to be non-empty is that (1.1.1) holds and 
$$\sum_{i=1}^k d_i \ \leq \ k(k-1) +\sum_{i=k+1}^n \min\left\{k, d_i \right\} \quad 
\text{for} \quad k=1, \ldots, n. \tag1.1.2$$

Our main goal is to estimate the cardinality $|G(D)|$ of $G(D)$. Using the obtained estimate, we 
deduce a concentration result for a random graph $G \in G(D)$ sampled from the uniform 
probability measure on $G(D)$. 

\subhead (1.2) The maximum entropy matrix and tame degree sequences
 \endsubhead The following matrix plays the crucial role 
in our construction.

Let us consider the space ${\Bbb R}^{n \choose 2}$ of vectors $x=\left(\xi_{\{j, k\}}\right)$, where 
$\{j, k\}$ is an unordered pair of indices $1 \leq j \ne k \leq n$. We consider the polytope 
${\Cal P} \subset {\Bbb R}^{n \choose 2}$, ${\Cal P}={\Cal P}(D)$, defined by the equations
$$\sum_{j: \ j \ne k} \xi_{\{j,k\}} =d_k \quad \text{for} \quad k=1,\ldots, n$$ 
and inequalities 
$$0 \leq  \xi_{\{j, k\}} \leq 1.$$
The integer points in ${\Cal P}(D)$ correspond to the labeled graphs with degree sequence $D$, 
which we write as 
$$G(D)={\Cal P}(D) \cap {\Bbb Z}^{n \choose 2}.$$

We assume that ${\Cal P}(D)$ is non-empty. 
We consider the following {\it entropy function} on ${\Cal P}(D)$:
$$\split H(x)=&\sum_{\{j, k\}}\left( \xi_{\{j, k\}} \ln {1 \over \xi_{\{j, k\}}} + \left(1-\xi_{\{j, k\}}\right) 
\ln {1 \over 1-\xi_{\{j, k\}}}\right) \\
&\text{for} \quad x=\left(\xi_{\{j, k\}} \right). \endsplit$$
Since $H$ is a strictly concave function, it attains its maximum on ${\Cal P}$ at a unique point,
$z=\left(\zeta_{\{j, k\}}\right)$, $z=z(D)$, which we call the 
{\it maximum entropy matrix} associated with the 
degree sequence $D$. Matrix $z$ can be easily calculated by interior point methods, see
\cite{NN94}. 

For $0< \delta \leq 1/2 $ we say that the degree sequence $D$ is $\delta$-{\it tame} if 
the polytope ${\Cal P}(D)$ is non-empty and if 
$$\delta \ \leq \ \zeta_{\{j ,k\}} \ \leq \ 1-\delta \quad \text{for all} \quad 1 \leq j \ne k \leq n,$$
where $z=\left(\zeta_{\{j, k\}}\right)$ is the maximum entropy matrix associated with degree
sequence $D$. In Theorem 2.1 we state some sufficient conditions for a degree sequence
$D$ to be tame. 

\subhead (1.3) Quadratic form $q$ and related quantities \endsubhead 
Let $z=\left(\zeta_{\{j, k\}}\right)$ be the maximum entropy matrix associated with a tame degree 
sequence $D$. We consider the following quadratic form $q: {\Bbb R}^n \longrightarrow {\Bbb R}$,
$$q(t)={1 \over 2} \sum_{\{j, k\}} \left(\zeta_{\{j, k\}}-\zeta_{\{j, k\}}^2 \right) \left(\tau_j +\tau_k\right)^2
\quad \text{for} \quad t=\left(\tau_1, \ldots, \tau_n \right).  \tag1.3.1$$
It is easy to see that $q$ is positive definite for $n >2$. Let us consider the Gaussian probability 
measure on ${\Bbb R}^n$ with density proportional to $e^{-q}$. We 
define the following random variables $f, h: {\Bbb R}^n \longrightarrow {\Bbb R}$, 
$$\aligned f(t)=&{1 \over 6} \sum_{\{j, k\}} \zeta_{\{j, k\}} \left(1 -\zeta_{\{j, k\}}\right) \left(2 \zeta_{\{j, k\}}-1 \right) \left(\tau_j + \tau_k\right)^3 \quad
\text{and} \\
h(t)=&{1 \over 24} \sum_{\{j, k\}} \zeta_{\{j, k\}} \left(1-\zeta_{\{j, k\}} \right) 
\left(6 \zeta_{\{j, k\}}^2 - 6 \zeta_{\{j, k\}} +1 \right) \left(\tau_j + \tau_k \right)^4 \\
& \qquad \qquad \text{for} \quad t=\left(\tau_1, \ldots, \tau_n \right). \endaligned \tag1.3.2$$
Let 
$$\mu=\EE f^2 \quad \text{and} \quad \nu=\EE h.$$

Our main result is as follows.
\proclaim{(1.4) Theorem} Let us fix $0 < \delta <1/2$. Let $D=\left(d_1, \ldots, d_n \right)$ be 
a $\delta$-tame degree sequence such that 
$d_1 + \ldots + d_n \equiv 0 \mod 2$, let $z=\left(\zeta_{\{j, k\}}\right)$ be the maximum 
entropy matrix as defined in Section 1.2 and let the quadratic form $q$ and values of $\mu$ and $\nu$ be as defined in Section 1.3. Let us define an $n \times n$ symmetric matrix $Q=\left(\omega_{jk}\right)$
by 
$$\split \omega_{jk}= &\zeta_{\{j, k\}} \left(1-\zeta_{\{j, k\}}\right) \quad \text{for} \quad j \ne k \quad \text{and} \\
\omega_{jj} =&d_j-\sum_{k: \ k \ne j} \zeta_{\{j, k\}}^2 \quad \text{for} \quad j=1, \ldots, n. \endsplit$$
Then $Q$ is positive definite and the value of
$$ {2 e^{H(z)} \over (2 \pi)^{n/2} \sqrt{\det Q}} \exp\left\{ -{\mu \over 2} +\nu\right\} \tag1.4.1$$
approximates the number of graphs $|G(D)|$ with degree sequence $D$ within a relative 
error which approaches 0 as $n \longrightarrow +\infty$. 

More precisely, for any $0<\epsilon \leq 1/2$ the value of (1.4.1) approximates $|G(D)|$ within 
relative error $\epsilon$ provided 
$$n \ \geq \ \left({1 \over \epsilon} \right)^{\gamma(\delta)},$$
where $\gamma=\gamma(\delta)$ is a positive constant.
\endproclaim 

The main term 
$${2 e^{H(z)} \over (2 \pi)^{n/2} \sqrt{\det Q}} \tag1.4.2$$ 
of formula (1.4.1) is the ``Gaussian approximation'' formula of \cite{BH10}, whose appearance, as 
is discussed in \cite{BH10}, is explained by the Local Central Limit Theorem, see also the 
discussion below. The factor 
$$\exp\left\{ -{\mu \over 2} +\nu \right\}$$ 
is the ``Edgeworth correction'' factor, see \cite{BH09b}. In the course of the proof of Theorem 1.4, we establish a two-sided bound 
$$\gamma_1(\delta) \ \leq \ \exp\left\{ -{\mu \over 2} +\nu \right\}\ \leq \gamma_2(\delta)$$
for some constants $\gamma_1(\delta), \gamma_2(\delta)>0$, as long as the 
degree sequence $D$ remains $\delta$-tame.

We note that computing the expectation of a polynomial with respect to the Gaussian probability 
measure is a linear algebra problem, cf. also Section 5.2. Hence apart from computing the maximum entropy 
matrix $z$, which can be done by interior point methods, computing the value of (1.4.1) is a
linear algebra problem which can be solved in $O(n^4)$ time in the unit cost model.

\subhead (1.5) Random graphs with prescribed degree sequences \endsubhead
Let us consider the set $G(D)$ of all labeled graphs with degree sequence $D$ as a 
finite probability space with the uniform measure. 
It is convenient to think of $G \in G(D)$ as of a subgraph of the complete graph $K_n$ with the 
set 
$$V=\bigl\{1, \ldots, n \bigr\}$$
of vertices and the set 
$$E=\Bigl\{ \{j, k\}: \quad 1 \leq j \ne k \leq n \Bigr\}$$
of edges.

Let us sample a graph $G \in G(D)$ at random.
What $G$ is likely to look like? 

As a corollary of Theorem 1.4, we prove that with overwhelming probability, for a random graph 
$G \in G(D)$ the number of edges of $G$ in a given set $S \subset E$
 with $|S|=\Omega(n^2)$  is very close to the sum of the entries of the maximum entropy matrix indexed by the elements of $S$.

\proclaim{(1.6) Theorem} Let us fix numbers $\kappa > 0$ and $0 < \delta  \leq 1/2$. Then 
there exists a number $\gamma(\kappa, \delta)>0$ such that the following holds. 

Suppose that $n \geq \gamma(\kappa, \delta)$ and that $D=\left(d_1, \ldots, d_n\right)$ is 
a $\delta$-tame degree sequence such that $d_1 + \ldots + d_n \equiv 0 \mod 2$.
For a set $S \subset E$, let $\sigma_S(G)$ be the 
number of edges of graph $G \in G(D)$ that belong to set $S$ and let 
$$\sigma_S(z)=\sum_{\{j, k\} \in S} \zeta_{\{j, k\}},$$
where $z=\left(\zeta_{\{j, k\}} \right)$ is the maximum entropy matrix.
Suppose that 
$|S| \ \geq \ \delta n^2$ and let 
$$\epsilon = \delta {\ln n \over \sqrt{n}}.$$
If $\epsilon \leq 1$ then for a uniformly chosen random graph $G \in G(D)$, we have 
$$\PP \Bigl\{G \in G(D): \quad (1-\epsilon) \sigma_S(z) \ \leq \ \sigma_S(G) \ \leq \ (1+\epsilon) 
\sigma_S(z) \Bigr\} \ \geq \ 1 - 2n^{-\kappa n}.$$
\endproclaim

The idea of the proof is as follows. For $1\leq j \ne k \leq n$, let $x_{\{j, k\}}$ be independent 
Bernoulli random variables such that 
$$\PP\bigl\{x_{\{j, k\}} =1 \bigr\} =\zeta_{\{j, k\}} \quad \text{and} \quad 
\PP\bigl\{x_{\{j, k\}} =0 \bigr\} =1-\zeta_{\{j, k\}}.$$
As is shown in \cite{BH10}, the probability mass function of the random vector $X=\left(x_{\{j, k\}}\right)$
is constant on the integer points of ${\Cal P}(D)$ and is equal to $e^{-H(z)}$ at each 
$G\in G(D)$, so that the vector $X$ conditioned on $G(D)$ is uniform. Theorem 1.4 then implies 
that the probability that $X \in G(D)$ is not too small. On the other hand, standard large 
deviation inequalities imply that the sum $\sum_{\{j, k\} \in S} x_{\{j, k\}}$ concentrates about the 
value of $\sigma_S(z)=\sum_{\{j, k\}} \zeta_{\{j, k\}}$. We supply the details of the proof in Section 10.

In many respects random graphs $G \in G(D)$ behave like random graphs on the set
$\{1, \ldots, n\}$ of vertices, with pairs $\{j, k\}$ chosen as the edges of $G$ independently with probabilities $\zeta_{\{j, k\}}$, where $z=\left(\zeta_{\{j, k\}}\right)$ is the maximum entropy matrix.
As is discussed in \cite{BH10}, the distribution of the multivariate Bernoulli 
random vector $X=\left(x_{\{j, k\}} \right)$ 
is the distribution of the largest entropy among all multivariate Bernoulli random vectors
constrained by 
$$\EE y_k = d_k \quad \text{for} \quad k=1, \ldots, n$$
where 
$$y_k =\sum_{j: \ j \ne k} x_{\{j, k\}}.$$

We remark that we obtain the ``Gaussian approximation'' term (1.4.2) if we assume
that the vector of random variables $Y=\left(y_1, \ldots, y_n \right)$
is asymptotically Gaussian around its expectation $\left(d_1, \ldots, d_n \right)$.
As it turns out, $Y$ is not exactly Gaussian but is not very far from it.

It looks plausible that both Theorem 1.4 and Theorem 1.6 can be extended to degree sequences 
$D$ allowing a moderate number of entries $\zeta_{\{j, k\}}$ of the maximum entropy matrix 
to be arbitrarily close to 1 or 0. Our proofs, however, do not seem to allow such an extension
with Theorem 4.1 being the main obstacle. Some of our proofs (mostly in Sections 5 and 6) are 
similar to those of \cite{BH09a}, where we applied the maximum entropy approach of \cite{BH10} 
to count non-negative integer matrices with prescribed row and column sums. 
\bigskip
The paper is organized as follows.

In Section 2, we give several examples and extensions concerning our main result, Theorem 1.4
and also discuss related work in the literature. 

In Section 3, we present an integral representation for the number $|G(D)|$ of graphs and 
also describe the plan of the proof of Theorem 1.4.

The rest of the paper deals with the proofs.

\head 2. Examples and extensions \endhead 

Sometimes one can tell that a degree sequence is tame without computing the maximum 
entropy matrix.
\proclaim{(2.1) Theorem} Let us fix real numbers $0 < \alpha < \beta < 1$ such that 
$$\beta < 2 \sqrt{\alpha} -\alpha, \quad \text{or, equivalently,} \quad (\alpha +\beta)^2 < 4\alpha.$$
Then there exists a real number $\delta=\delta(\alpha, \beta) >0$ and 
 a positive integer $n_0=n_0(\alpha, \beta)$ such that any 
degree sequence $D=\left(d_1, \ldots, d_n \right)$ satisfying 
$$\alpha \ < \ {d_i \over n-1} \ < \ \beta \quad \text{for} \quad i=1, \ldots, n$$ 
is $\delta$-tame provided $n > n_0$.

One can choose 
$$\split &n_0= \max\left\{ {\beta \over \alpha(1-\beta)}, \quad {4(\beta-\alpha) \over 4 \alpha -(\alpha +\beta)^2} \right\} +1
\quad \text{and} \\
&\delta={\epsilon^6 \over 1+\epsilon^6} \quad \text{where} \quad 
\epsilon =\min\left\{ \alpha, \quad \alpha - {(\alpha+\beta)^2 \over 4} \right\}. \endsplit$$
\endproclaim 

For example, degree sequences $D=\left(d_1, \ldots, d_n \right)$ satisfying 
$$0.25\ < \ {d_i\over n-1} \ < \ 0.74  \quad \text{for} \quad i=1, \ldots, n$$
or 
$$0.01\ < \ {d_i \over n-1} \ < \ 0.18 \quad \text{for} \quad i=1, \ldots, n$$ 
or 
$$0.81 \ < \ {d_i \over n-1} \ < \ 0.89 \quad \text{for} \quad i=1, \ldots, n$$
are $\delta$-tame for some $\delta >0$ and all sufficiently large $n$.

We prove Theorem 2.1 in Section 12. 

\example{(2.2) On the boundary of $\delta$-tameness}
 Let us choose rational $0 < \alpha < \beta < 1$ such that 
$$\beta = 2 \sqrt{\alpha} -\alpha. \tag2.2.1$$
Clearly, $\beta > \alpha$. Let us choose a positive integer $n$ such that 
$\alpha n$ and $\beta n$ are even integers and let us consider the degree sequence 
$$d_1 = \ldots = d_k =\beta n \quad \text{and} \quad d_{k+1} = \ldots = d_n =\alpha n 
\quad \text{for} \quad k=n\sqrt{\alpha}$$
(note that $k$ is necessarily integral).
The Erd\H{o}s-Gallai condition (1.1.2) for $k=n\sqrt{\alpha}$, reduces to 
$$\beta \ \leq \ 2 \sqrt{\alpha} - \alpha - {1 \over n}. \tag2.2.2$$
In particular, (2.2.1) does not even guarantee that the polytope ${\Cal P}(D)$ is non-empty.

In \cite{JSM92} Jerrum, Sinclair and McKay discuss under what conditions
an approximation formula for $|G(D)|$ which depends ``smoothly'' on $D$ may exist.
They
 describe the phenomenon of the number of 
graphs $|G(D)|$ changing sharply when the degree sequence $D$ is varying only slightly 
around some special values of $D$. 
This phenomenon is apparently explained by the fact that the dimension of the polytope ${\Cal P}(D)$ may change abruptly or the polytope may 
disappear altogether when $D$ lies on the boundary of the Erd\H{o}s-Gallai conditions (1.1.2).
Theorem 8.1 of \cite{JSM92} states that for
$$d_+=\max\bigl\{d_i,\ i =1, \ldots, n \bigr\} \quad \text{and} \quad d_-=\min\bigl\{d_i,\ i=1, \ldots, n \bigr\},$$
as long as 
$$\left(d_+ - d_- +1\right)^2 \ \leq \ 4 d_-\left(n-d_+ -1\right), \tag2.2.3$$
the degree sequence $D$ is $P$-{\it stable}, meaning that increasing one of the degrees $d_i$
and decreasing another by 1 does not change $|G(D)|$ by more than a factor of $n^{10}$
(this, in turn, implies that there are polynomial time randomized approximation algorithms for 
computing $|G(D)|$ and sampling a random graph $G \in G(D)$). The condition of our Theorem 
2.1 is only marginally stronger than (2.2.3).
 
As Sourav Chatterjee pointed out to us, 
Lemma 4.1 of recent \cite{CDS11} shows that a sequence $D$ is $\delta$-tame provided it lies 
sufficiently deep inside the polyhedron defined by 
the Erd\H{o}s-Gallai conditions (1.1.2). 
 
Our example shows that the bounds of Theorem 2.1 are essentially the best possible if 
we take into account only the largest and the smallest degree of a vertex of the graph.
\endexample
 
\subhead (2.3) Regular graphs \endsubhead
In \cite{MW90} McKay and Wormald compute the asymptotic of $|G(D)|$ for regular graphs,
where
$$d_1= \ldots = d_n =d,$$
and almost regular graphs, where 
$$\left|d_i -d \right| < n^{{1\over 2} +\epsilon} \quad \text{for} \quad i=1, \ldots, n$$
for a sufficiently small $\epsilon >0$; see also \cite{McK11} for recent developments and \cite{McK10} for a survey.
One can show that the formula of Theorem 1.4 is equivalent to the asymptotic formula of 
\cite{MW90} for regular or almost regular graphs.

In the case of regular graphs, symmetry requires that 
$$\zeta_{\{j, k\}} = {d \over n-1} \quad \text{for all} \quad 1 \leq j \ne k \leq n$$
for the maximum entropy matrix $z=\left(\zeta_{\{j,k\}}\right)$.

\subhead (2.4) Approximations in the cut norm \endsubhead
The {\it cut norm} (sometimes called the {\it normalized} cut norm) of a real $m \times n$ matrix $A=\left(a_{jk}\right)$ is defined by 
$$\|A\|_{\text{cut}}={1 \over mn} \max_{J, K} \left| \sum \Sb j \in J, \ k \in K \endSb a_{jk} \right|,$$
where the maximum is taken over all non-empty subsets $J \subset \{1,\ldots, m\}$ and
$K \subset \{1, \ldots, n\}$. Let us choose set $S$ in Theorem 1.6 of the form 
$$S=\Bigl\{\{j, k\}: \quad j \in J, \ k \in K, \ j\ne k \Bigr\}
\qquad \text{for some} \qquad J, K \subset \{1, \ldots, n\}.$$ 
We note that there 
are not more than $2^{2n}$ distinct sets $S$ of this form. Theorem 1.6 implies
 that as $n$ grows, the maximum 
entropy matrix $z(D)$ approximates the adjacency matrix of the overwhelming majority of graphs
 $G \in G(D)$ within an 
error of $O\left(n^{-1/2} \ln n \right)$ in the cut norm.

Shortly after the first version of this paper appeared, using a different approach, Chatterjee, Diaconis and Sly \cite{CDS11}
described {\it graph limits} of graphs from $G(D)$ as $n$ grows. A graph limit is a certain 
function on $[0, 1] \times [0, 1]$, viewed as an ``infinite matrix'', which naturally arises as a limit object for a Cauchy sequence 
in the cut norm of adjacency matrices of graphs \cite{LS06}.  
Graph limits constructed in \cite{CDS11} can indeed be viewed as ``infinite maximum entropy
matrices''.

\example{(2.5) Enumeration of bipartite graphs} A natural version of the problem concerns enumeration
of labeled {\it bipartite} graphs with a given degree sequence or, equivalently, $m \times n$ matrices 
with 0-1 entries and prescribed row sums $R=\left(r_1, \ldots, r_m\right)$ and column sums 
$C=\left(c_1, \ldots, c_n\right)$. 
We assume that 
$$r_1 + \ldots + r_m =c_1 +\ldots + c_n.$$
A simple necessary and sufficient condition for a 0-1 matrix with prescribed row and column sums 
to exist is given by the Gale-Ryser Theorem, see, for example, Corollary 6.2.5 of \cite{BR91}.

Let us consider the polytope ${\Cal P}(R, C)$ of $m \times n$ matrices $x=\left(\xi_{jk}\right)$ defined by 
the equations 
$$\sum_{k=1}^n \xi_{jk}=r_j \quad \text{for} \quad j=1, \ldots, m \qquad \text{and} \qquad
\sum_{j=1}^m \xi_{jk}=c_k \quad \text{for} \quad k=1, \ldots, n$$
and inequalities
$$0 \ \leq \ \xi_{jk} \ \leq \ 1 \quad \text{for all} \quad j, k.$$
Let us compute the {\it maximum entropy} matrix $z=\left(\zeta_{jk}\right)$ as 
the necessarily unique matrix $z \in {\Cal P}(R, C)$ that maximizes
$$H(x)=\sum_{jk} \left( \xi_{jk} \ln {1 \over \xi_{jk}} +\left(1-\xi_{jk}\right) \ln {1 \over 1-\xi_{jk}} \right) 
\quad \text{for} \quad x =\left(\xi_{jk}\right)$$ 
on ${\Cal P}(R,C)$.
For $0 < \delta \leq 1/2$, we say that the margins $(R, C)$ are $\delta$-{\it tame} if 
$$\delta m \ \leq \ n \quad \text{and} \quad \delta n \ \leq \ m$$ 
and 
$$\delta \ \leq \ \xi_{jk} \ \leq \ 1-\delta \quad \text{for all} \quad j, k.$$
Suppose that the margins $(R, C)$ are indeed $\delta$-tame for some $\delta>0$.
Let us define a quadratic form $q: {\Bbb R}^{m+n} \longrightarrow {\Bbb R}$ 
by 
$$\split q(s, t)=&{1 \over 2} \sum_{j, k} \left( \zeta_{jk} -\zeta_{jk}^2\right) \left(\sigma_j +\tau_k \right)^2 
\\ & \text{for} \quad (s, t)=\left(\sigma_1, \ldots, \sigma_m; \tau_1, \ldots, \tau_n \right). 
\endsplit \tag2.5.1$$
Let
$$u=\left(\underbrace{1, \ldots, 1}_{\text{$m$ times}}; \underbrace{-1, \ldots, -1}_{\text{$n$ times}}\right)
\tag2.5.2$$
and let $L=u^{\bot}$ be the orthogonal complement to $u$ in ${\Bbb R}^{m+n}$. 
Then the restriction $q|L$ of $q$ onto $L$
is strictly positive definite and we define $\det q|L$ as the product of the non-zero eigenvalues 
of $q$. We consider the Gaussian probability measure on $L$ with density proportional to 
$e^{-q}$ and define random variables $f,g: L \longrightarrow {\Bbb R}$ by 
$$\split f(s, t)=&{1 \over 6} \sum_{j,k} \zeta_{jk} \left(1 -\zeta_{jk}\right) \left(2 \zeta_{jk}-1 \right) \left(\sigma_j + \tau_k\right)^3 \quad
\text{and} \\
h(s, t)=&{1 \over 24} \sum_{j,k} \zeta_{jk} \left(1-\zeta_{jk} \right) 
\left(6 \zeta_{jk}^2 - 6 \zeta_{jk} +1 \right) \left(\sigma_j + \tau_k \right)^4 \\
& \qquad \qquad \text{for} \quad (s, t)=\left(\sigma_1, \ldots, \sigma_m; \tau_1, \ldots, \tau_n \right). \endsplit \tag2.5.3$$
We define
$$\mu=\EE f^2 \quad \text{and} \quad \nu=\EE h.$$
Then the number $|R, C|$ of 0-1 matrices with row sums $R$ and column sums $C$ is 
$$|R, C| ={e^{H(z)} \sqrt{m+n} \over (4 \pi)^{(m+n-1)/2} \sqrt{\det q|L}} \exp\left\{-{\mu \over 2} + \nu \right\}
\Bigl(1+o(1)\Bigr) \tag2.5.4$$
provided $m, n \longrightarrow +\infty$ in such a way that the margins $(R, C)$ remain $\delta$-tame 
for some $\delta>0$. We sketch the proof of (2.5.4) in Section 11.

Canfield and McKay \cite{CM05} obtained an asymptotic formula of $|R, C|$ when all row 
sums are equal, $r_1 =\ldots = r_m$ and all column sums are equal, $c_1=\ldots =c_n$, which was 
later extended to the case of ``almost equal'' row sums and ``almost equal'' column sums 
\cite{CGM08}, see also \cite{GM09}. The maximum entropy matrix $z$ was introduced in 
\cite{Ba10} where a cruder asymptotic formula 
$$\ln |R, C| \approx H(z)$$
was established without the $\delta$-tameness assumption and for a wider class of enumeration 
problems, including enumeration of 0-1 matrices with prescribed row and column sums 
{\it and} zeros in prescribed position. It was also shown in \cite{Ba10} that a random matrix 0-1 
with prescribed row and column sums concentrates about the maximum entropy matrix $z$.
\endexample

\head 3. An integral representation for the number of graphs \endhead

In \cite{BH10} we proved the following general result; see Theorem 5, Lemma 11 and formula (16)
there.

\proclaim{(3.1) Theorem} Let $P \subset {\Bbb R}^p$ be a polyhedron defined by the system 
of linear equations $Ax=b$, where $A$ is a $n \times p$ matrix with columns 
$a_1, \ldots, a_p \in {\Bbb Z}^n$ and $b \in {\Bbb Z}^n$ is an integer vector, and inequalities
$0 \leq x \leq 1$ (the inequalities are understood coordinate-wise). Suppose that $P$ 
has a non-empty interior, that is, contains a point $x=\left(\xi_1, \ldots, \xi_p\right)$ such that 
$0 < \xi_j < 1$ for $j=1, \ldots, p$.

Then the function 
$$H(x)=\sum_{j=1}^p \left( \xi_j \ln {1 \over \xi_j}  +\left(1-\xi_j \right) \ln {1 \over 1-\xi_j} \right)
\quad \text{for} \quad x=\left(\xi_1, \ldots, \xi_p \right)$$
attains its maximum on $P$ at a unique point $z=\left(\zeta_1, \ldots, \zeta_p\right)$ such that 
$0 <\zeta_j <1$ for $j=1, \ldots, p$. 

Let us consider the parallelepiped $\Pi=[-\pi, \pi]^n$, $\Pi \subset {\Bbb R}^n$. 
Then the number $|P \cap \{0, 1\}^p|$ of 0-1 points in $P$ can be written as 
$$| P \cap \{0, 1\}^p |={e^{H(z)} \over (2 \pi)^n} \int_{\Pi} e^{-i \langle t, b \rangle} 
\prod_{j=1}^p \left(1 -\zeta_j + \zeta_j e^{i \langle a_j, t \rangle} \right) \ dt,$$
where $\langle \cdot, \cdot \rangle$ is the standard scalar product in ${\Bbb R}^n$, 
$dt$ is the standard Lebesgue measure in ${\Bbb R}^n$ and $i=\sqrt{-1}$. 
\endproclaim 
{\hfill \hfill \hfill} \qed 

The idea of the proof is as follows. Let $X=\left(x_1, \ldots, x_p \right)$ be a random vector 
of independent Bernoulli random variables such that $\PP\left\{x_j=1\right\}=\zeta_j$
and $\PP\left\{x_j=0\right\}=1-\zeta_j$ for $j=1, \ldots, p$. It turns out that the probability mass function 
of $X$ is constant on the set $P \cap \{0, 1\}^p$ and equals $e^{-H(z)}$ for every 0-1 point in
$P$. Letting $Y=AX$, we obtain 
$$|P \cap \{0, 1\}^p| =e^{H(z)} \PP \{X \in P \} =e^{H(z)} \PP\{Y=b\}$$
and the probability in question is written as the integral of the characteristic function of $Y$.

Since 
$$\sum_{j=1}^p \zeta_j a_j =b,$$
in a neighborhood of the origin $t=0$ the integrand can be written as 
$$\aligned e^{-i \langle t, b \rangle} 
&\prod_{j=1}^p \left(1 -\zeta_j + \zeta_j e^{i \langle a_j, t \rangle} \right) =
\\ \exp\Biggl\{&-{1 \over 2}\sum_{j=1}^p  \zeta_j \left(1-\zeta_j \right) \langle a_j, t \rangle^2 \\
                        &+{i \over 6} \sum_{j=1}^p \zeta_j \left(1-\zeta_j \right) \left(2 \zeta_j -1 \right) 
                        \langle a_j, t \rangle^3 \\
                        &+{1 \over 24} \sum_{j=1}^p \zeta_j \left(1-\zeta_j \right)\left(6 \zeta_j^2 -6 \zeta_j +1\right)
                        \langle a_j, t \rangle^4 \\
                        &+O\left( \sum_{j=1}^p \left(\zeta_j +1 \right)^5 \langle a_j, t \rangle^5\right) \Biggr\}.
  \endaligned \tag3.2$$
Note that the linear term is absent in the expansion.

We obtain the following corollary.

\proclaim{(3.3) Corollary} Let $D=\left(d_1, \ldots, d_n \right)$ be a degree sequence such that 
the polytope ${\Cal P}(D)$ defined in Section 1.2 has a non-empty interior and let $z=\left(\zeta_{\{j, k\}} \right)$ be 
the maximum entropy matrix. 
Let 
$$\split F(t)=&\exp\left\{-i \sum_{m=1}^n d_m \tau_m \right\}
\prod_{\{j, k\}} \left( 1-\zeta_{\{j, k\}} +\zeta_{\{j, k\}} e^{i \left(\tau_j + \tau_k \right)}  \right) \\
&\quad \text{for} \quad t=\left(\tau_1, \ldots, \tau_n \right). \endsplit$$
Then for the parallelepiped $\Pi=[-\pi, \pi]^n$, we have
$$|G(D)| = {e^{H(z)} \over (2 \pi)^n} \int_{\Pi} F(t) \ dt.$$
\endproclaim
\demo{Proof} Follows by Theorem 3.1.
{\hfill \hfill \hfill} \qed
\enddemo

We note that in the case of regular and almost regular graphs (see Section 2.3) the integral 
of Corollary 3.3 is the same as the one evaluated by McKay and Wormald \cite{MW90}.

\subhead (3.4) Plan of the proof of Theorem 1.4 \endsubhead
We use the integral representation of Corollary 3.3. Let us define subsets
${\Cal U}, {\Cal W} \subset \Pi$ by 
$${\Cal U}=\Bigl\{ \left(\tau_1, \ldots, \tau_n\right): \quad \left| \tau_j\right| \ \leq \ {\ln n \over \sqrt{n}} 
\quad \text{for} \quad j =1, \ldots, n \Bigr\}$$
and 
$$\split {\Cal W}=&\Bigl\{ \left(\tau_1, \ldots, \tau_n\right): \quad \left| \tau_j -\sigma_j \pi \right| \ \leq \ {\ln n \over \sqrt{n}} \quad \text{for some} \quad  \sigma_j=\pm 1 \\ &\qquad \qquad  \text{and} \quad j =1, \ldots, n \Bigr\}.
\endsplit$$
We show that the integral of $F(t)$ over $\Pi \setminus ({\Cal U} \cup {\Cal W})$ is asymptotically negligible.
Namely, in Section 8 we prove that the integral
$$\int_{\Pi \setminus ({\Cal U} \cup {\Cal W})} |F(t)| \ dt $$ is asymptotically negligible compared to the integral
$$\int_{{\Cal U}} |F(t)| \ dt. \tag3.4.1$$
It is easy to show that 
$$\int_{\Cal U} F(t) \ dt =\int_{\Cal W} F(t) \ dt,$$
provided $d_1 + \ldots + d_n$ is even.

In Section 7,  we evaluate 
$$\int_{\Cal U} F(t) \ dt. \tag3.4.2$$
In particular, we show that the integrals (3.4.1) and (3.4.2) have the same order of magnitude and so the 
integral of $F(t)$ outside of ${\Cal U} \cup {\Cal W}$ is indeed asymptotically irrelevant.

From (3.2) one can deduce that asymptotically as $n \longrightarrow +\infty$,
$$F(t) \approx \exp\bigl\{ -q(t) + i f(t) + h(t) \bigr\}  \quad \text{for} \quad t \in {\Cal U},$$
where $q$ is defined by (1.3.1) and $f$ and $h$ are defined by (1.3.2).

Let us consider the Gaussian probability measure in ${\Bbb R}^n$ with density proportional 
to $e^{-q}$. In Section 6, we prove that with respect to that measure 
$$h(t) \approx \EE h =\nu \quad \text{almost everywhere in} \quad {\Cal U}. \tag3.4.3$$ 
This allows us to conclude that 
$$\int_{\Cal U}  \exp\bigl\{ -q(t) + i f(t) + h(t) \bigr\} \ d t \approx e^{\nu} \int_{\Cal U} \exp\bigl\{-q(t) + i f(t) \bigr\} \ dt.$$
In Section 5, we prove that asymptotically, as $n \longrightarrow +\infty$, function $f$ is a Gaussian 
random variable, so 
$$\aligned \int_{\Cal U} \exp\bigl\{-q(t) + i f(t)\bigr\} \ dt \approx &\int_{{\Bbb R}^n} \exp\bigl\{-q(t) + i f(t) \bigr\} \ dt \\
\approx &\exp\left\{ -{1 \over 2} \EE f^2 \right\}  \int_{{\Bbb R}^n} e^{-q(t)} \ dt, \endaligned 
\tag3.4.4$$
which concludes the evaluation of (3.4.2).

The crucial consideration used in proving (3.4.3) and (3.4.4) is that with respect to the Gaussian 
probability measure in ${\Bbb R}^n$ with density proportional to $e^{-q}$, the coordinate functions 
$\tau_1, \ldots, \tau_n$ are weakly correlated, that is,
$$\split \left| \EE \tau_j \tau_k \right| =&O\left({1 \over n^2}\right) \quad \text{for} \quad j \ne k \quad \text{and} \\
\EE \tau_j^2 = &O\left({1 \over n}\right) \quad \text{for} \quad j=1, \ldots, n. \endsplit \tag3.4.5$$
We prove (3.4.5) in Section 4, where we essentially use the $\delta$-tameness assumption.

\subhead (3.5) Notation \endsubhead By $\gamma$, sometimes with an index or a list of parameters,
we denote a positive constant depending only on the listed parameters. The most common 
appearance will be $\gamma(\delta)$, a positive constant depending only on the 
$\delta$-tameness constant $\delta$.

As usual, for two functions $g_1$ and $g_2$, where $g_2$ is non-negative, we write 
$g_1=O(g_2)$ if $|g_1|  \leq \gamma g_2$ and 
$g_1 =\Omega(g_2)$ if $g_1 \geq \gamma g_2$ for some $\gamma>0$.

\head 4. Correlations \endhead 

Let $z=\left(\zeta_{\{j, k\}} \right)$ be the maximum entropy matrix as defined in Section 1.2. 
We assume that 
$$0 \ < \ \zeta_{\{j, k\}} \ < \ 1 \quad \text{for all} \quad j \ne k.$$
We define the quadratic form $q: {\Bbb R}^n \longrightarrow {\Bbb R}$ 
by 
$$q(t)={1 \over 2} \sum_{\{j, k\}} \left(\zeta_{\{j, k\}} -\zeta_{\{j, k\}}^2 \right) \left(\tau_j + \tau_k\right)^2
\quad \text{for} \quad t=\left(\tau_1, \ldots, \tau_n\right).$$
For $n>2$ the  quadratic form $q$ is strictly positive definite. We consider the Gaussian
probability measure on ${\Bbb R}^n$ with density proportional to $e^{-q}$.
We consider a point $t=\left(\tau_1, \ldots, \tau_n \right)$ as a random vector and 
$\tau_1, \ldots, \tau_n$ as random variables. 

The main result of this section is as follows.
\proclaim{(4.1) Theorem} For any $0< \delta \leq 1/2$ there exists $\gamma(\delta)>0$ such that the 
following holds.

Suppose that 
$$\delta \ \leq \ \zeta_{\{j, k\}}\ \leq \ 1-\delta \quad \text{for all} \quad j \ne k.$$
Then 
$$ \split \left| \EE \tau_j \tau_k \right| \ &\leq \ {\gamma(\delta) \over n^2} \quad \text{provided} \quad 
j \ne k \quad \text{and} \\
\EE \tau_j^2  \ &\leq \ {\gamma(\delta) \over n} \quad \text{for} \quad j=1, \ldots, n. \endsplit$$
\endproclaim

We will often consider the following situation. Let $\psi: {\Bbb R}^n \longrightarrow {\Bbb R}$ 
be a positive definite quadratic form. We consider the Gaussian probability measure 
in ${\Bbb R}^n$ with density proportional to $e^{-\psi}$. For a polynomial (random variable) 
$f: {\Bbb R}^n \longrightarrow {\Bbb R}$ we denote by $\EE(f; \ \psi)$ its expectation with respect 
to the measure. For a subspace $L \subset {\Bbb R}^n$, we consider the restriction $\psi|L$ of 
$\psi$ onto $L$ and the Gaussian probability measure on $L$ with density proportional to 
$e^{-\psi|L}$. For a polynomial $f: {\Bbb R}^n \longrightarrow {\Bbb R}$, we denote by $\EE(f; \ \psi|L)$ the expectation 
of the restriction $f: L \longrightarrow {\Bbb R}$ with respect to that  Gaussian probability measure on 
$L$. We will use the following standard fact: 
suppose that ${\Bbb R}^n = L_1 \oplus L_2$ is a decomposition of ${\Bbb R}^n$ into the direct 
sum of orthogonal subspaces such that 
$$\psi\left(t_1 + t_2\right) = \psi\left(t_1\right) + \psi\left(t_2\right) \quad 
\text{for all} \quad t_1 \in L_1 \quad \text{and} \quad t_2 \in L_2,$$
so that the coordinates $t_1 \in L_1$  and $t_2 \in L_2$ of the point $t=t_1+t_2$, $t \in {\Bbb R}^n$
are independent.
Let $\ell_1, \ell_2: {\Bbb R}^n \longrightarrow {\Bbb R}$ be linear functions. 
Then 
$$\EE \left( \ell_1 \ell_2; \ \psi\right)= \EE\left(\ell_1 \ell_2; \ \psi|L_1 \right) + 
\EE\left(\ell_1 \ell_2; \ \psi|L_2 \right).$$
Indeed, writing $t=t_1+t_2$ with $t_1 \in L_1$ and $t_2 \in L_2$ and noting that 
$\ell_{1,2}(t)=\ell_{1, 2}(t_1) +\ell_{1,2}(t_2)$,
we obtain
$$\split \EE \left(\ell_1(t) \ell_2(t); \ \psi \right)=&
\EE\left(\ell_1(t_1) \ell_2(t_1); \ \psi \right) + \EE\left( \ell_1(t_1) \ell_2(t_2);\  \psi \right) 
\\ &\quad + \EE\left( \ell_1(t_2) \ell_2(t_1); \ \psi \right) +\EE \left( \ell_1(t_2) \ell_2(t_2); \ \psi \right) \\=
&\EE\left( \ell_1 \ell_2; \ \psi| L_1 \right) + 2\EE \left(\ell_1; \ \psi\right) \EE\left(\ell_2; \ \psi \right) +
\EE \left(\ell_1 \ell_2; \ \psi| L_2 \right)\\ =
& \EE\left(\ell_1 \ell_2; \ \psi| L_1 \right) + \EE \left(\ell_1 \ell_2; \ \psi| L_2 \right).\endsplit $$

We deduce Theorem 4.1 from the following result.

\proclaim{(4.2) Proposition} Let $n >2$ and let $\xi_{\{j, k\}}$, $1 \leq j \ne k \leq n$ be a set 
of numbers such that
$$\alpha \ \leq \ \xi_{\{j ,k\}} \ \leq \beta \quad \text{for all} \quad j, k$$
and some $\beta > \alpha>0$.

Let 
$$\sigma_k =\sum_{j: \ j \ne k} \xi_{\{j, k\}} \quad \text{for} \quad k=1, \ldots, n.$$
Let us consider the quadratic form $\psi: {\Bbb R}^n  \longrightarrow {\Bbb R}$
defined by 
$$\psi(t) ={1 \over 2} \sum_{\{j, k\}} \xi_{\{j ,k\}}\left( {\tau_j \over \sqrt{\sigma_j}} + {\tau_k \over \sqrt{\sigma_k}} \right)^2
\quad \text{for} \quad t=\left(\tau_1, \ldots, \tau_n \right),$$
where the sum is taken over all unordered pairs of indices $1 \leq j \ne k \leq n$.
Then $\psi$ is a positive definite quadratic form and we consider the Gaussian probability 
measure in ${\Bbb R}^n$ with density proportional to $e^{-\psi}$. 

Let 
$$\epsilon={\alpha \over \beta}.$$
Then for $n > 2/\epsilon$ we have 
$$\split &\left| \EE \tau_j \tau_k \right|  \ \leq \ {n^2 \over \epsilon^{5/2} (n-\epsilon)(n \epsilon -2) (n-1)}
+{3 \over 2 \epsilon n} \quad \text{provided} \quad j \ne k  \quad \text{and} \\
 &\left| \EE \tau_j^2 -1  \right| \leq \ {n^2 \over \epsilon^{5/2} (n-\epsilon)(n \epsilon -2) (n-1)}
+{3 \over 2 \epsilon n} \quad \text{for} \quad j=1, \ldots, n. \endsplit$$
\endproclaim
\demo{Proof} Clearly, $\psi$ is positive definite. Let
$$v=\left(\sqrt{\sigma_1}, \ldots, \sqrt{\sigma_n} \right).$$
Then $v$ is an eigenvector of $\psi$ with eigenvalue 1.
Indeed, the gradient of $\psi$ at $t=v$ is $2v$:
$${\partial \over \partial \tau_j} \psi(t) \Big|_{t=v}= {1 \over \sqrt{\sigma_j}} \sum_{k: \ k \ne j} 2 \xi_{\{j, k\}} =
2 \sqrt{\sigma_j} \quad \text{for} \quad j=1, \ldots, n.$$
Let 
$$L=v^{\bot} \subset {\Bbb R}^n$$ 
be the orthogonal complement to $v$.
Hence $L$ is defined in ${\Bbb R}^n=\bigl\{\left(\tau_1,  \ldots, \tau_n \right) \bigr\}$ by the equation 
$$\sum_{j=1}^n \tau_j \sqrt{\sigma_j} =0.$$
We write 
$$\psi(t)={1 \over 2} \sum_{j=1}^n \tau_j^2 + \sum_{\{j, k\}} {\xi_{\{j, k\}} \over \sqrt{\sigma_j \sigma_k}} \tau_j \tau_k.$$
We remark that 
$$\alpha (n-1) \ \leq \ \sigma_j \ \leq \beta(n-1) \quad \text{for} \quad j=1, \ldots, n. \tag4.2.1$$
Let 
$$\omega=\sum_{j=1}^n \sigma_j \ \geq \ \alpha n(n-1)$$
and let us define the quadratic form $\phi: {\Bbb R}^n \longrightarrow {\Bbb R}$ by 
$$\phi(t) ={1 \over  \omega} \left( \sum_{j=1}^n \tau_j \sqrt{\sigma_j}  \right)^2 \quad \text{for} \quad 
t=\left(\tau_1, \ldots, \tau_n \right).$$
Hence $\phi(t)$ is a form of rank 1 and $v$ is an eigenvector of $\phi$ with eigenvalue 1.

We define a perturbation 
$$\tilde{\psi} = \psi -{\epsilon^2 \over 2} \phi \quad \text{for} \quad \epsilon={\alpha \over \beta}.$$
Hence $\tilde{\psi}$ is a positive definite quadratic form such that
$$\tilde{\psi}(t) =\psi(t) \quad \text{for all} \quad t \in L,$$
and $v$ is an eigenvector of $\tilde{\psi}$ with eigenvalue $1-\epsilon^2/2$.

Let us consider the Gaussian probability measure on ${\Bbb R}^n$ with density 
proportional to $e^{-\tilde{\psi}}$. Our immediate goal is to estimate the covariances  
 $\EE \tau_j \tau_k$ with respect to that measure.

Denoting by $\langle \cdot, \cdot \rangle$ the standard scalar product in ${\Bbb R}^n$, we
can write 
$$\tilde{\psi}(t)={1 \over 2} \langle (I +Q) t,\ t \rangle,$$
where $I$ is the $n \times n$ identity matrix and
$Q=\left(q_{jk} \right)$ is an $n \times n$ symmetric matrix such that $v$ is an eigenvector of $Q$ 
with eigenvalue $1-\epsilon^2$. We have 
$$\split q_{jk} =&{\xi_{\{j, k\}} \over \sqrt{\sigma_j \sigma_k}} - \epsilon^2  {\sqrt{\sigma_j \sigma_k} \over \omega} \quad 
\text{for} \quad j \ne k \quad \text{and} \\
q_{jj}=&-{\epsilon^2 \sigma_j \over \omega} \quad \text{for} \quad j=1, \ldots, n. \endsplit$$

It follows by (4.2.1) that 
$$\aligned {1 \over \epsilon(n-1)}  \ \geq \ &q_{jk} \ \geq  0 \quad \text{for} \quad j \ne  k \quad \text{and} \\
0 \ \geq \ &q_{jj} \ \geq \ -{\epsilon \over n} \quad \text{for} \quad j=1, \ldots, n. \endaligned $$
The covariance matrix $R=\left( \EE \tau_j \tau_k; \ \tilde{\psi} \right)$ of the Gaussian measure with 
density proportional to $e^{-\tilde{\psi}}$ is 
$$\split (I + Q)^{-1} =&\left(\left( 1-{\epsilon \over n} \right)I+ \left( {\epsilon \over n} I + Q \right) \right)^{-1} 
  =\left(1 -{\epsilon \over n} \right)^{-1} (I + P)^{-1}, \quad \text{where} \quad \\
  &P=\left(1-{\epsilon \over n} \right)^{-1} \left({\epsilon \over n} I +Q \right). \endsplit$$
Hence $P=\left(p_{jk}\right)$ is a symmetric matrix such that 
$$0 \ \leq \ p_{jk} \ \leq \ \left(1-{\epsilon \over n} \right)^{-1} {1 \over \epsilon(n-1)} \quad \text{for all} \quad j, k. \tag4.2.2$$
Furthermore, $v$ is an eigenvector of $P$ with eigenvalue $(1-\epsilon^2+\epsilon/n)/(1-\epsilon/n)$, 
so
$$P v =\lambda v \quad \text{for} \quad \lambda=
\left(1-{\epsilon \over n} \right)^{-1} \left(1 -\epsilon^2 +{\epsilon \over n} \right). \tag4.2.3$$
Let us bound the entries of a positive integer power 
$P^d=\left(p_{jk}^{(d)} \right)$ of $P$.
Let 
$$\kappa = \left(1-{\alpha \over \beta n} \right)^{-1} {\beta \over \alpha^{3/2} (n-1)^{3/2}} \quad \text{and let} \quad y =\kappa v, \quad y=\left(\eta_1, \ldots, \eta_n \right).$$
By (4.2.1) and (4.2.2), we have 
$$p_{j k} \ \leq \eta_j \quad \text{for all} \quad j, k. \tag4.2.4$$
Also, by (4.2.1) we have 
$$\eta_j \ \leq \ \left(1-{\epsilon \over n} \right)^{-1} {1 \over \epsilon^{3/2} (n-1)} \quad \text{for all} \quad j.\tag4.2.5$$
Furthermore, $y$ is an eigenvector of $P$ with eigenvalue $\lambda$ defined by (4.2.3),
and hence $y$ is an eigenvector of the $d$-th power $P^d=\left(p_{jk}^{(d)}\right)$ with eigenvalue $\lambda^d$.
Combining this with (4.2.4) and (4.2.5), for $d \geq 0$ we obtain
$$p_{jk}^{(d+1)}=\sum_{m=1}^n p_{jm}^{(d)} p_{mk} \ \leq \ \sum_{m=1}^n p_{jm}^{(d)} \eta_m 
=\lambda^d \eta_j \ \leq \ \lambda^d \left(1-{\epsilon \over n} \right)^{-1} {1 \over \epsilon^{3/2} (n-1)}.$$
We note that for $n > 2/\epsilon$ we have $0 <\lambda <1$. Consequently, the series 
$$(I+P)^{-1} =I +\sum_{d=1}^{+\infty} (-1)^d P^d $$
converges absolutely and we can bound the entries of the matrix 
$$R=(I+Q)^{-1}=\left(1 -{\epsilon \over n} \right)^{-1}(I +P)^{-1},$$
$R=\left(r_{jk} \right)$ by 
$$\split & \left|r_{jk}\right| \ \leq \ \left(1 -{\epsilon \over n} \right)^{-2} {1 \over \epsilon^{3/2}(n-1)} {1 \over 1-\lambda}
\quad \text{if} \quad j \ne k \quad \text{and} \\
&\left| r_{jj} -1\right|  \ \leq \  \left(1 -{\epsilon \over n} \right)^{-2} {1 \over \epsilon^{3/2}(n-1)} {1 \over 1-\lambda} \quad \text{for} \quad j=1, \ldots, n. \endsplit$$
We have 
$${1 \over 1-\lambda} =\left(1-{\epsilon \over n} \right) \left(\epsilon^2 -{2 \epsilon \over n} \right)^{-1}.$$
Since $R$ is the covariance matrix of the Gaussian probability measure with density 
proportional to $e^{-\tilde{\psi}}$ , we obtain
$$\aligned &\left| \EE \left(\tau_j \tau_k; \ \tilde{\psi}\right)  \right| \ \leq \ {n^2 \over \epsilon^{5/2}  (n-\epsilon) (n \epsilon -2) (n-1)}
\quad \text{provided} \quad j \ne k \quad \text{and} \\
&\left| \EE \left(\tau_j^2; \ \tilde{\psi}\right) -1 \right| \ \leq \ {n^2 \over \epsilon^{5/2}  (n-\epsilon) (n \epsilon -2) (n-1)} 
\quad \text{for} \quad j=1, \ldots, n. \endaligned \tag4.2.6$$

Now we go back to the form $\psi$ and the Gaussian probability measure with 
density proportional to $e^{-\psi}$. Since $v$ is an eigenvector of both $\psi$ and $\tilde{\psi}$, since $L=u^{\bot}$ and since $\psi$ and $\tilde{\psi}$ coincide on $L$, for
any linear functions $\ell_1, \ell_2: {\Bbb R}^n \longrightarrow {\Bbb R}$, we have 
$$\aligned \EE\left(\ell_1 \ell_2; \ \psi \right)=&\EE\left( \ell_1 \ell_2;\ \psi|L \right) + 
\EE \left( \ell_1 \ell_2; \ \psi| \spa(v) \right) \\
=&\EE\left(\ell_1 \ell_2; \ \tilde{\psi}| L \right) + \EE \left( \ell_1 \ell_2; \ \psi| \spa(v) \right) \\
=&\EE\left(\ell_1 \ell_2; \ \tilde{\psi} \right) -\EE\left(\ell_1 \ell_2; \ \tilde{\psi}| \spa(v) \right)\\ &\qquad +
\EE \left(\ell_1 \ell_2; \ \psi|\spa(v) \right). \endaligned \tag4.2.7$$
We note that the gradient of the coordinate function $\tau_j$ restricted to $\spa(v)$ is 
$\sqrt{\sigma_j/\omega}$. Since $v$ is an 
eigenvector of $\psi$ with eigenvalue 1 and an eigenvector of $\tilde{\psi}$ with eigenvalue 
$1-\epsilon^2/2$, we have
$$\split &\EE\left(\tau_j \tau_k; \ \psi| \spa(v)\right) ={\sqrt{\sigma_j \sigma_k} \over 2 \omega} \quad \text{and} \\
&\EE \left( \tau_j \tau_k; \ \tilde{\psi}|\spa(v) \right) = {\sqrt{\sigma_j \sigma_k} \over (2-\epsilon^2)\omega}  \quad \text{for all} \quad j, k. \endsplit$$
By (4.2.1) we have 
$$\left| \EE\left(\tau_j \tau_k; \ \psi| \spa(v)\right) \right| \ \leq \  {1 \over 2 \epsilon n} \quad \text{and} \quad
\left| \EE \left( \tau_j \tau_k; \ \tilde{\psi}|\spa(v) \right) \right| \ \leq \  {1 \over \epsilon n}.$$
The proof now follows by (4.2.6) and (4.2.7).
{\hfill \hfill \hfill} \qed
\enddemo 

Now we are ready to prove Theorem 4.1.

\demo{Proof of Theorem 4.1} 
Let us define 
$$\xi_{\{j, k\}}=\zeta_{\{j, k\}}-\zeta^2_{\{j, k\}} \quad \text{for all} \quad j \ne k$$
and let us choose $\alpha=\delta-\delta^2$ and $\beta=1/4$ in Proposition 4.2.
We define $\sigma_j$ and $\psi$ as in Proposition 4.2
and consider a linear transformation
$$\left(\tau_1, \ldots, \tau_n\right) \longmapsto \left(\tau_1 \sqrt{\sigma_1}, \ldots, \tau_n \sqrt{\sigma_n} \right).$$
Then the push-forward of the Gaussian probability measure with density proportional to 
$e^{-q}$ is the Gaussian probability measure with density proportional to $e^{-\psi}$.
Therefore,
$$\EE \left( \tau_j \tau_k; \ q\right) = {1 \over \sqrt{\sigma_j \sigma_k}}  \EE\left(\tau_j \tau_k; \ \psi \right).$$
Since 
$$\sigma_j \ \geq \ \delta \alpha (n-1) \quad \text{for} \quad j=1, \ldots, n,$$
The proof follows by Proposition 4.2. 
{\hfill \hfill \hfill} \qed
\enddemo 

We will need the following lemma.

\proclaim{(4.3) Lemma}
 Let $q_0: {\Bbb R}^n \longrightarrow {\Bbb R}$, $n \geq 2$, be the quadratic form 
defined by the formula 
$$q_0(t)={1 \over 2} \sum_{\{j, k\}} \left(\tau_j + \tau_k \right)^2 \quad 
\text{for} \quad t=\left(\tau_1, \ldots, \tau_n\right).$$
Then the eigenspaces of $q_0$ are as follows:
the 1-dimensional eigenspace $E_1$ with eigenvalue $n-1$ spanned by the vector 
$u=(1, \ldots, 1)$ and 
the $(n-1)$-dimensional eigenspace $E_2=u^{\bot}$  with eigenvalue $(n-2)/2$.
\endproclaim
\demo{Proof} We have
$${\partial \over \partial \tau_k} q(t) \big|_{t=u}=2n-2.$$
Hence the gradient of $q_0(t)$ at $t=u$ is $(2n-2)u$, so $u$ is an eigenvector with 
eigenvalue $(n-1)$. For $ t \in u^{\bot}$ we have $\tau_1 + \ldots +\tau_n=0$ and 
hence 
$${\partial \over \partial \tau_k} q(t) = \sum_{j:\ j \ne k} \left(\tau_j +\tau_k \right)=(n-2) \tau_k.$$
Therefore, the gradient of $q_0(t)$ at $t \in u^{\bot}$ is $(n-2)t$, and so $t$ is an eigenvector 
with eigenvalue $(n-2)/2$.
{\hfill \hfill \hfill} \qed
\enddemo

\head 5. The third degree term \endhead

The main result of this section is the following theorem.

\proclaim{(5.1) Theorem} For unordered pairs $\{j, k\}$, $1 \leq j \ne k \leq n$, let 
$u_{\{j, k\}}$ be Gaussian random variables such that 
$$\EE u_{\{j, k\}}=0 \quad \text{for all} \quad j, k.$$ 
Suppose further that for some $\theta >0$ we have 
$$\EE u_{\{j, k\}}^2 \ \leq \ {\theta \over n} \quad \text{for all} \quad j, k$$
and that 
$$\left| \EE u_{\{j_1, k_1\}} u_{\{j_2, k_2\}} \right| \ \leq \ {\theta \over n^2} \quad 
\text{provided} \quad \{j_1, k_1\} \cap \{j_2, k_2\} =\emptyset.$$
Let 
$$U=\sum_{\{j, k\}} u_{\{j, k\}}^3.$$
Then for some constant $\gamma(\theta)>0$ and any $0< \epsilon <1/2$ we have
$$\left| \EE \exp\left\{i U\right\} -\exp\left\{-{1 \over 2} \EE U^2 \right\}\right| \ \leq \ \epsilon$$
provided
$$n \ \geq \ \left({1 \over \epsilon} \right)^{\gamma(\theta)}.$$
Furthermore,
$$\EE U^2 \ \leq \ \gamma(\theta)$$ 
for some $\gamma(\theta)>0$. Here $i=\sqrt{-1}$.
\endproclaim

We apply Theorem 5.1 in the following situation. 
Let $q: {\Bbb R}^n \longrightarrow {\Bbb R}$ be the quadratic form defined by (1.3.1).

Let us consider the Gaussian probability measure on ${\Bbb R}^n$ with density proportional 
to $e^{-q}$. We define random variables $u_{\{j, k\}}$ by 
$$u_{\{j, k\}}(t) =\root 3 \of{{1 \over 6} \zeta_{\{j,k\}}\left(1-\zeta_{\{j, k\}}\right) \left(2 \zeta_{\{j, k\}}-1\right)} 
\left(\tau_j + \tau_k\right) \quad \text{for} \quad t=\left(\tau_1, \ldots, \tau_n\right).$$
Then for the function $f(t)$ defined by (1.3.2) we have 
$$f=\sum_{\{j, k\}} u_{\{j, k\}}^3.$$

In this section, all implied constants in the ``$O$'' notation are absolute.

 Our main tool is Wick's formula 
for the expectation of the product of random Gaussian variables.  

\subhead (5.2) Wick's formula \endsubhead Let $w_1, \ldots, w_l$ be Gaussian random variables
such that 
$$\EE w_1 = \ldots =\EE w_l=0.$$
Then 
$$\split \EE \left(w_1 \cdots w_l \right)=&0 \quad \text{if} \quad l \quad \text{is odd} \quad \text{and} \\
\EE\left(w_1 \ldots w_l \right)=&\sum \left(\EE w_{i_1} w_{i_2} \right) \cdots \left(\EE w_{i_{l-1}} w_{i_l}
\right) \quad \text{if} \quad l=2r \quad \text{is even}, \endsplit$$
where the sum is taken over all $(2r)!/r! 2^r$ unordered partitions of the set of indices 
$\{1, \ldots, l\}$ into $r=l/2$ pairwise disjoint unordered pairs 
$\left\{i_1, i_2\right\}, \ldots, \left\{i_{l-1}, i_l \right\}$, see for example, \cite{Zv97}.
 Such a partition is called a {\it matching} of 
the random variables $w_1, \ldots, w_l$ and we say that $w_i$ and $w_j$ are {\it matched} if 
they form a pair in the matching. 

In particular,
$$\EE w^{2r} = {(2r)! \over r! 2^r} \left( \EE w^2 \right)^r \tag5.2.1$$ 
for a centered Gaussian random variable $w$. We will also use 
that 
$$\EE w_1^3 w_2^3 =9 \left(\EE w_1^2 \right) \left(\EE w_2^2 \right) \left(\EE w_1 w_2\right) +
6 \left( \EE w_1 w_2 \right)^3 \tag5.2.2$$ and later in Section 6 that 
$$\aligned \cov\left(w_1^4, w_2^4 \right)=&\EE \left(w_1^4 w_2^4\right)-\left(\EE w_1^4 \right)
\left(\EE w_2^4\right) \\ =
&72 \left( \EE w_1 w_2 \right)^2 \left( \EE w_1^2 \right)\left(\EE w_2^2\right)+
24 \left(\EE w_1 w_2 \right)^4. \endaligned \tag5.2.3$$

\subhead (5.3) Representing monomials by graphs \endsubhead
Let $x_{\{j, k\}}: \quad 1 \leq j \ne k \leq n$
be formal commuting variables. We interpret a monomial in $x_{\{j, k\}}$ as a weighted graph as follows.
Let $K_n$ be the complete graph with vertices $1, \ldots, n$ and edges $\{j, k\}$ for 
$1 \leq j \ne k \leq n$. A {\it weighted graph} $G$ is a set of edges $\{j, k\}$ of $K_n$ with positive 
integer weights $\alpha_{\{j, k\}}$ on them. The set of vertices of $G$ consists of all vertices of the edges of $G$.  
With $G$, we associate a monomial 
$$m_G(x)=\prod_{\{j, k\} \in G} x_{\{j, k\}}^{\alpha_{\{j, k\}}}.$$
The {\it weight} of $G$ is the degree of $m_G(x)$, that is, $\sum_{\{j, k\} \in G} \alpha_{\{j, k\}}$. 
We observe that for any $p$ there are not more than $r^{O(r)} n^p$ distinct weighted graphs $G$ of weight $2r$ 
on $p$ vertices. 

In what follows, given a set of random variables, we construct auxiliary Gaussian random variables with the same 
matrix of covariances. This is always possible since the matrix of covariances is positive semi-definite.

Our proof of Theorem 5.1 is based on the following combinatorial lemma.

\proclaim{(5.4) Lemma} For the Gaussian random variables $u_{\{j, k\}}$ of Theorem 5.1, 
let us introduce auxiliary Gaussian random variables $v_{\{j, k\}}$ such that 
$$\split \EE v_{\{j, k\}} = &0 \quad \text{for all} \quad 1\leq j \ne k \leq n \qquad \text{and} \\
\EE v_{\{j_1, k_1\}} v_{\{j_2, k_2\}} = &\EE u_{\{j_1, k_1\}}^3 u_{\{j_2, k_2\}}^3 \quad 
\text{for all} \quad 1 \leq j_1 \ne k_1, \ j_2 \ne k_2 \leq n. \endsplit$$
Given a weighted graph $G$ of weight $2r$, $r>1$, let us represent it as a vertex-disjoint 
union 
$$G=G_0 \cup G_1,$$
where $G_0$ consists of $s$ isolated edges of weight 1 each and $G_1$ is a graph with 
no isolated edges of weight 1 (we may have $s=0$ and $G_0$ empty).

Then
\roster 
\item We have 
$$\left| \EE m_G\left(u_{\{j, k\}}^3: \ 1 \leq j \ne k \leq n\right) \right|, 
\  \left| \EE m_G \left(v_{\{j, k\}}: \ 1 \leq j \ne k \leq n\right) \right| \ \leq \ {r^{O(r)} \theta^{3r} \over n^{3r+s/2}}$$
if $s$ is even and 
$$\left| \EE m_G\left(u_{\{j, k\}}^3: \ 1 \leq j \ne k \leq n \right)\right|, 
\  \left| \EE m_G \left(v_{\{j, k\}}: \ 1 \leq j \ne k \leq n \right)\right| \ \leq \ {r^{O(r)} \theta^{3r} \over n^{3r+(s+1)/2}}$$
if $s$ is odd.
\item  If $s$ is even and $G_1$ is a vertex-disjoint union of $r-s/2$ connected components, each 
consisting of a pair of edges of weight 1 sharing exactly one common vertex, then 
$$\left| \EE m_G\left(u_{\{j, k\}}^3: \ 1 \leq j \ne k \leq n \right) - \EE m_G\left(v_{\{j, k\}}: \ 1 \leq j \ne k \leq n\right)\right| \ \leq \ {r^{O(r)} \theta^{3r} 
\over n^{3r + s/2 +1}}.$$
\item 
If $s$ is even and $G_1$ is a vertex-disjoint union of $r-s/2$ connected components, each consisting
of a pair of edges of weight 1 sharing exactly one common vertex, then $G$ has 
precisely $3r + s/2$ vertices. In all other cases, $G$ has strictly fewer than $3r+s/2$ vertices.
\endroster
\endproclaim 
\demo{Proof} If $\{j_1, k_1\} \cap \{j_2, k_2\} =\emptyset$ we say that the pair of variables 
$u_{\{j_1, k_1\}}$, $u_{\{j_2, k_2\}}$ and the pair of variables $v_{\{j_1, k_1\}}$,
$v_{\{j_2, k_2\}}$ are {\it weakly correlated}. If $\{j_1, k_1\} \cap \{j_2, k_2\} \ne \emptyset$
we say that the pairs of variables are {\it strongly correlated}. Pairs of variables indexed by edges in 
different connected components of $G$ are necessarily weakly correlated. 

To prove Part (1) we use Wick's formula of Section 5.2. By (5.2.2), we obtain
$$\aligned &\EE \left(v_{\{j_1, k_1\}} v_{\{j_2, k_2\}}\right) = O\left({\theta^3 \over n^4}\right)
\quad \text{if the pair} \quad v_{\{j_1, k_1\}}, \ v_{\{j_2, k_2\}} \\ &\qquad \qquad \text{is weakly correlated}, \\
&\EE \left(v_{\{j_1, k_1\}} v_{\{j_2, k_2\}}\right) = O\left({\theta^3 \over n^3} \right) 
\quad \text{if the pair} \quad v_{\{j_1, k_1\}}, \ v_{\{j_2, k_2\}} \\ &\qquad \qquad
\text{is strongly correlated.}  \endaligned \tag5.4.1$$
Since for each isolated edge $\{j_1, k_1\} \in G_0$ variable $v_{\{j_1, k_1\}}$
 has to be matched with variable $v_{\{j_2, k_2\}}$ indexed by an edge $\{j_2, k_2\}$ in a 
 different connected component, we conclude that every matching of the 
 set 
 $$\Bigl\{ v_{\{j, k\}}: \quad \{j, k\} \in G \Bigr\} \tag5.4.2$$
 contains at least $s/2$ weakly correlated pairs and hence
 $$\left| \EE m_G \left(v_{\{j, k\}}: \ 1 \leq j \ne k \leq n\right) \right| \ \leq \ r^{O(r)} 
 \left( {\theta^3 \over n^4} \right)^{s/2} \left( {\theta^3 \over n^3} \right)^{r-s/2}.$$
 Moreover, if $s$ is odd, then the number of weakly correlated pairs is at least $(s+1)/2$
 and hence 
  $$\left| \EE m_G \left(v_{\{j, k\}}: \ 1 \leq j \ne k \leq n\right) \right| \ \leq \ r^{O(r)} 
 \left( {\theta^3 \over n^4} \right)^{(s+1)/2} \left( {\theta^3 \over n^3} \right)^{r-(s+1)/2}.$$
 Similarly, since for each isolated edge $\{j_1, k_1\} \in G_0$ at least one copy of the 
 variable $u_{\{j_1, k_1\}}$ has to be matched with a copy of variable $u_{\{j_2, k_2\}}$ 
 indexed by an edge in a different connected component, we conclude that 
 every matching of the multiset 
 $$\Bigl\{ u_{\{j, k\}}, u_{\{j, k\}}, u_{\{j, k\}}: \quad \{j, k\} \in G \Bigr\} \tag5.4.3$$
 contains at least $s/2$ weakly correlated pairs, and hence
 $$\left| \EE m_G \left(u_{\{j, k\}}^3: \ 1 \leq j \ne k \leq n\right) \right| \ \leq \ r^{O(r)} 
 \left( {\theta \over n^2} \right)^{s/2} \left( {\theta \over n} \right)^{3r-s/2}.$$
  Moreover, if $s$ is odd, then the number of weakly correlated pairs is at least $(s+1)/2$
 and hence 
  $$\left| \EE m_G \left(u_{\{j, k\}}^3: \ 1 \leq j \ne k \leq n\right) \right| \ \leq \ r^{O(r)} 
 \left( {\theta \over n^2} \right)^{(s+1)/2} \left( {\theta \over n} \right)^{3r-(s+1)/2}.$$
 This concludes the proof of Part (1).
 
 To prove Part (2), let us define $\Sigma_v(G)$ as the sum in the Wick's formula for 
 $\EE m_G\left(v_{\{j, k\}}: \ 1 \leq j \ne k \leq n\right)$ taken over all matchings of the set (5.4.2) of the following structure:
 we split the edges of $G$ into $r$ pairs, pairing each isolated edge with another isolated 
 edge and pairing each edge in a connected component of $G$ consisting of two edges with the remaining edge 
 in the same connected component. Then we match every variable $v_{\{j_1, k_1\}}$ with the 
 variable $v_{\{j_2, k_2\}}$ such that $\{j_1, k_1\}$ is paired with $\{j_2, k_2\}$.
 Reasoning as in the proof of Part (1), we conclude that 
 $$\left| \EE m_G \left(v_{\{j, k\}}: \ 1 \leq j \ne k \leq n\right) - \Sigma_v(G) \right| \ \leq \ 
 {r^{O(r)} \theta^{3r} \over n^{3r + s/2+1}},$$
 since every matching of the set (5.4.2) which is not included in $\Sigma_v(G)$ contains 
 at least $s/2+1$ weakly correlated pairs.
 
 Similarly, let us define $\Sigma_u(G)$ as the sum in the Wick's formula for \break
 $\EE m_G\left(u_{\{j, k\}}^3: \ 1 \leq j \ne k \leq n \right)$ taken over all matchings of the multiset (5.4.3) of the 
 following structure: we split the edges of $G$ into $r$ pairs as above and match every 
 copy of variable $u_{\{j_1, k_1\}}$ with a copy of variable $u_{\{j_2, k_2\}}$ indexed by an 
 edge in the same pair (in particular, we may match copies of the same variable). Reasoning as 
 in the proof of Part (1), we conclude that 
 $$\left| \EE m_G \left(u_{\{j, k\}}^3: \ 1 \leq j \ne k \leq n\right) - \Sigma_u(G) \right| \ \leq \ 
 {r^{O(r)} \theta^{3r} \over n^{3r + s/2+1}},$$
 since every matching of the multiset (5.4.3) which is not included in $\Sigma_u(G)$ 
 contains at least $s/2+1$ weakly correlated pairs.
 
 The proof of Part (2) follows since 
 $$\Sigma_u(G)=\Sigma_v(G)$$ by Wick's formula.
 
To prove Part (3), we note that a connected weighted graph $G$ of weight $e$ contains a spanning 
tree with at most $e$ edges and hence has at most $e+1$ vertices. In particular, a connected 
graph $G$ of weight $e$ contains fewer than $3e/2$ vertices unless $G$ is an 
isolated edge of weight $1$ or a pair of edges of weight 1 each, sharing one common vertex.
Therefore, $G$ has at most 
$$2s+ {3 \over 2} (2r-s)=3r + {s \over 2}$$
vertices and strictly fewer vertices, unless $s$ is even and the connected components 
of $G_1$ are pairs of edges of weight 1 each sharing one common vertex.
{\hfill \hfill \hfill} \qed
\enddemo

\subhead (5.5) Proof of Theorem 5.1 \endsubhead Let 
$v_{\{j, k\}}$ be Gaussian random variables defined in Lemma 5.4 and let 
$$V=\sum_{\{j, k\}} v_{\{j, k\}}.$$
Since there are $O\left(n^3\right)$ strongly correlated pairs $v_{j_1, k_1}, v_{j_2, k_2}$ and
there are $O\left(n^4\right)$ weakly correlated pairs, by (5.4.1) we have 
$$\EE V^2 = \EE U^2 =O\left(\theta^3\right). \tag5.5.1$$
Since $V$ is a Gaussian random variable, we have 
$$\EE e^{iV} =\exp\left\{ - {1 \over 2} \EE V^2 \right\} =\exp\left\{ -{1 \over 2} \EE U^2 \right\}.
\tag5.5.2$$
Our goal is to show that  $\EE e^{iV}$ and $\EE e^{iU}$ are asymptotically equal as 
$n \longrightarrow +\infty$. 

By symmetry, the odd moments of $U$ and $V$ are 0:
$$\EE U^k = \EE V^k =0 \quad \text{if} \quad k >0 \quad \text{is odd}. \tag5.5.3$$
The even moments of $U$ and $V$ can be expressed as 
$$\aligned \EE U^{2r}=&\sum_G a_G \EE m_G \left(u_{\{j, k\}}^3:\ 1\leq j \ne k \leq n \right) \\ 
\EE V^{2r}=&\sum_G a_G \EE m_G  \left(v_{\{j, k\}}:\ 1 \leq j \ne k \leq n \right), \endaligned$$
where the sum is taken over all weighted graphs $G$ of weight $2r$ and 
$$ 1 \ \leq \ a_G \ \leq (2r)!.$$
Let ${\Cal G}_{2r}$ be the set of weighted graphs $G$ of weight $2r$ whose connected components 
are an even number $s$ of isolated edges of weight 1 and $r-s/2$ pairs of edges of weight 1 
sharing one common vertex. 
Since there are no more than $r^{O(r)}n^p$ distinct weighted graphs of weight $2r$ with $p$
vertices, by Parts (1) and (3) of Lemma 5.4, we have
$$\split \left| \EE U^{2r} -\sum_{G \in {\Cal G}_{2r}} a_G \EE m_G \left(u_{\{j, k\}}^3:\ 1 \leq j \ne k \leq n \right) \right| 
\ \leq \ &{r^{O(r)} \theta^{3r} \over n} \quad \text{and} \\
\left| \EE V^{2r} -\sum_{G \in {\Cal G}_{2r}} a_G \EE m_G \left(v_{\{j, k\}}:\ 1 \leq j \ne k \leq n \right) \right| 
\ \leq \ &{r^{O(r)} \theta^{3r} \over n}. \endsplit$$
Therefore, by Part (2) of Lemma 5.4, 
$$\left| \EE U^{2r} - \EE V^{2r} \right| \ \leq \  {r^{O(r)} \theta^{3r}\over n}. \tag5.5.4$$
From Taylor's Theorem 
$$\left| e^{ix} -\sum_{s=0}^{2r-1} i^s {x^s \over s!} \right| \ \leq \ {x^{2r} \over (2r)!} \quad \text{for}
\quad x \in {\Bbb R},$$ it follows that
$$\left| \EE e^{iU} - \EE e^{iV} \right| \ \leq \ {\EE U^{2r} \over (2r)!} + {\EE V^{2r} \over (2r)!} 
+\sum_{s=0}^{2r-1} {|\EE U^s - \EE V^s| \over s!}.$$
 From (5.5.1) and (5.2.1) we deduce 
that for a positive integer $r$ we have 
$$\EE V^{2r}  \ \leq \ {(2r)! 2^{O(r)} \theta^{3r} \over r!}.$$
Therefore, by (5.5.4)
$$\left| \EE e^{iU} -\EE e^{iV}\right| \ \leq \ {2^{O(r)} \theta^{3r} \over r!} + {r^{O(r)} \theta^{3r} \over n}.
\tag5.5.5$$
Given $0 \leq \epsilon \leq 1/2$, one can choose a positive integer $r$ such that 
$$r \ln r = O\left(\theta^2 \ln {1 \over \epsilon} \right)$$
so that the first term in the right hand side of (5.5.5) does not exceed $\epsilon/2$.
It follows then that for all 
$$n \ \geq \ \left({1 \over \epsilon} \right)^{\gamma(\theta)}$$ we have 
$$\left| \EE e^{iU} - \EE e^{iV} \right| \ \leq \ \epsilon,$$
and the proof follows by (5.5.2).
{\hfill \hfill \hfill} \qed

\head 6. The fourth degree term \endhead

The main result of this section is the following theorem.

\proclaim{(6.1) Theorem} For unordered pairs $\{j, k\}$, $1 \leq j \ne k \leq n$, let $w_{\{j, k\}}$ 
be Gaussian random variables such that 
$$\EE w_{\{j, k\}}=0 \quad \text{for all} \quad j, k,$$
and let $\sigma_{\{j, k\}} \in \{-1, 1\}$ be numbers.

Suppose further that for some $\theta>0$ we have 
$$\EE w_{\{j, k\}}^2 \ \leq \ {\theta \over n} \quad \text{for all} \quad j, k,$$
and that 
$$\left| \EE w_{\{j_1, k_1\}} w_{\{j_2, k_2\}}  \right| \ \leq \ {\theta \over n^2} \quad \text{provided} 
\quad \{j_1, k_1\} \cap \{j_2, k_2\} =\emptyset.$$

Let 
$$W=\sum_{\{j, k\}} \sigma_{\{j, k\}} w_{\{j, k\}}^4.$$
Then for some constant $\gamma(\theta)>0 $ we have 
\roster
\item 
$$\EE |W|  \ \leq \ \gamma(\theta);$$
\item 
$$\var W \ \leq \ {\gamma(\theta) \over n};$$
\item 
$$\PP\Bigl\{|W| \ > \ \gamma(\theta) \Bigr\} \ \leq \ \exp\left\{-n^{1/5}\right\}$$
provided $n \geq \gamma_1(\theta)$ for some constant $\gamma_1(\theta)>0$.
\endroster
\endproclaim
We apply Theorem 6.1 in the following situation. Let $q: {\Bbb R}^n \longrightarrow {\Bbb R}$ 
be the quadratic form defined by (1.3.1).
Let us consider the Gaussian probability measure on ${\Bbb R}^n$ with density proportional 
to $e^{-q}$. We define random variables $w_{\{j, k\}}$ by 
$$\split w_{\{j, k\}}(t) = &\root 4 \of{{1 \over 24}\zeta_{\{j, k\}} \left(1-\zeta_{\{j, k\}}\right) \left| 6 \zeta_{\{j, k\}}^2 - 6 \zeta_{\{j, k\}} +1 \right|}
\left(\tau_j +\tau_k\right) \\ &\qquad \qquad  \text{for} \quad t=\left(\tau_1, \ldots, \tau_n\right)
\endsplit$$
and let 
$$\sigma_{\{j, k\}}=\sgn \left( 6 \zeta_{\{j, k\}}^2 - 6 \zeta_{\{j, k\}} +1\right).$$
Then for the function $h$ defined by (1.3.2),  we have 
$$h =\sum_{\{j, k\}} \sigma_{\{j, k\}} w_{\{j, k\}}^4.$$

While the proof of Parts (1)--(2) is done by a straightforward computation, to prove Part (3) we need
reverse H\"older inequalities for polynomials with respect to the Gaussian measure. 
\proclaim{(6.2) Lemma} Let $p$ be a polynomial of degree $d$ in random Gaussian variables 
$w_1, \ldots, w_l$. Then for $r>2$ we have 
$$\left(\EE |p|^r \right)^{1/r} \ \leq \ r^{d/2} \left( \EE p^2 \right)^{1/2}.$$
\endproclaim
\demo{Proof} This is Corollary 5 of \cite{Du87}. 
{\hfill \hfill \hfill} \qed
\enddemo

\subhead (6.3) Proof of Theorem 6.1 \endsubhead
All implied constants in the ``$O$'' notation below are absolute. 

By formula (5.2.1), 
$$\EE w_{\{j, k\}}^4 = 3 \left(\EE w_{\{j, k\}}^2\right)^2 = O\left({\theta^2 \over n^2}\right) $$ and hence 
$$\EE |W|=O\left(\theta^2\right)$$
and Part (1) follows.
Furthermore,
$$\var W = \sum \Sb \{j_1, k_1\} \\ \{j_2, k_2\} \endSb  \sigma_{\{j_1, k_1\}} \sigma_{\{j_2, k_2\}} 
\cov \left( w_{\{j_1, k_1\}}^4, w_{\{j_2, k_2\}}^4 \right).$$
By (5.2.3)  we have 
$$\cov \left(w_{\{j_1, k_1\}}^4,  w_{\{j_2, k_2\}}^4 \right) \ = \ O\left({\theta^4 \over n^4} \right)$$ 
and, additionally,
$$\cov \left(w_{\{j_1, k_1\}}^4,  w_{\{j_2, k_2\}}^4 \right) \ = \ O\left({\theta^4 \over n^6} \right)
\quad \text{provided} \quad \{j_1, k_1\} \cap \{j_2, k_2\} =\emptyset.$$ 
Therefore,
$$\var W = O\left({\theta^4 \over n}\right),  \tag6.3.1$$
which proves Part (2).

Finally, applying Lemma 6.2 with $d=4$ we deduce from (6.3.1) that for any $r>2$ 
$$\EE \left| W-\EE W \right|^r \ \leq \ r^{2r} n^{-r/2} 2^{O(r)} \theta^{2r}.$$
Choosing $r=n^{1/5}$, we conclude that for all sufficiently large $n \geq n_0(\theta)$
we have
$$\EE \left| W-\EE W \right|^r \ \leq \ \exp\left\{ -n^{1/5} \right\}.$$
By Markov's inequality, we obtain 
$$\PP\Bigl\{ |W -\EE W| >1  \Bigr\} \ \leq \ \exp\left\{-n^{1/5}\right\}$$
for all sufficiently large $n$ and the proof follows from Part (1).
{\hfill \hfill \hfill} \qed

\head 7. Computing the integral over a neighborhood of the origin \endhead 

We consider the integral 
$$\int_{\Pi} F(t) \ dt$$
of Corollary 3.3. Hence 
$$\split F(t)=&\exp\left\{-i \sum_{m=1}^n d_m \tau_m \right\} \prod_{\{j, k\}} 
\left(1-\zeta_{\{j, k\}} + \zeta_{\{j, k\}} e^{i(\tau_j + \tau_k)} \right) \\ &\text{for} \quad
t=\left(\tau_1, \ldots, \tau_n \right), \endsplit$$
where $D=\left(d_1, \ldots, d_n\right)$ is a given degree sequence, $z=\left(\zeta_{\{j, k\}}\right)$
is the maximum entropy matrix and $\Pi$ is the parallelepiped $[-\pi, \pi]^n$.
 We recall that the quadratic form $q: {\Bbb R}^n \longrightarrow {\Bbb R}$ 
is defined by 
$$q(t)={1 \over 2} \sum_{\{j, k\}} \left( \zeta_{\{j, k\}} - \zeta_{\{j, k\}}^2 \right) \left(\tau_j +\tau_k \right)^2
\quad \text{for} \quad t=\left(\tau_1, \ldots, \tau_n\right).$$
In this section, we prove the following main result.
\proclaim{(7.1) Theorem} Let us fix a number $0 < \delta \leq 1/2$ and suppose that 
$$\delta \ \leq \ \zeta_{\{j, k\}} \ \leq \ 1-\delta \quad \text{for all} \quad j \ne k.$$
Let $f, h: {\Bbb R}^n \longrightarrow {\Bbb R}$ be polynomials defined by (1.3.2). 
Let us define a neighborhood of the origin ${\Cal U} \subset \Pi$ by 
$${\Cal U}=\left\{ \left(\tau_1, \ldots, \tau_n\right): \quad \left| \tau_k\right| \ \leq \ {\ln n \over \sqrt{n}} 
\quad \text{for} \quad k=1, \ldots, n \right\}.$$
Let 
$$\Xi=\int_{{\Bbb R}^n} e^{-q(t)} \ dt $$ 
and let us consider the Gaussian probability measure in ${\Bbb R}^n$ with density 
$\Xi^{-1} e^{-q}$. Let 
$$\mu = \EE f^2 \quad \text{and} \quad \nu =\EE h.$$
Then 
\roster
\item 
$$\Xi \ \geq \ \left({4 \pi \over n}\right)^{n/2};$$
\item 
$$\mu, \ |\nu|, \ \EE |h|  \ \leq \ \gamma(\delta)$$
for some constant $\gamma(\delta)>0$;
\item For any $0<\epsilon \leq 1/2$
$$\left|  \int_{\Cal U} |F(t)| \ dt  - \exp\{\nu\} \  \Xi \right| 
 \ \leq \ \epsilon\  \Xi $$
provided 
$$n \ \geq \ \left({1 \over \epsilon}\right)^{\gamma(\delta)}$$
for some constant $\gamma(\delta)>0$;
\item For any $0 <\epsilon \leq 1/2$
$$\left| \int_{\Cal U} F(t) \ dt - \exp\left\{-{\mu \over 2} +\nu \right\} \ \Xi \right| 
 \ \leq \ \epsilon  \ \Xi$$
provided
$$n \ \geq \ \left( {1 \over \epsilon} \right)^{\gamma(\delta)}$$
for some $\gamma(\delta)>0$.
\endroster
\endproclaim
\demo{Proof}

In what follows, all constants implied in the ``$O$'' and ``$\Omega$'' notation depend only on the parameter 
$\delta$.

Let 
$$q_0(t)={1 \over 2} \sum_{\{j, k\}} \left(\tau_j + \tau_k \right)^2 \quad \text{for} \quad 
t=\left(\tau_1, \ldots, \tau_n \right)$$ 
as in Lemma 4.3. 
Then $q(t) \ \leq \ {1 \over 4} q_0(t)$ and hence 
$$\int_{{\Bbb R}^n} e^{-q} \ dt \ \geq \ \int_{{\Bbb R}^n} e^{-{1 \over 4}q_0(t)} \ dt =
\pi^{n/2}  \sqrt{4 \over n-1} \left({8 \over n-2}\right)^{n-1 \over 2} \ \geq \ \left({4 \pi  \over n}\right)^{n/2},$$
which proves Part (1).

Let us think of the coordinate functions $\tau_j$ as random variables with respect to 
the Gaussian probability measure with density proportional to $e^{-q}$. 

By Theorem 4.1, we have 
$$\aligned \left| \EE \tau_j \tau_k \right| \ = \ &O\left({1 \over n^2}\right) \quad \text{provided} \quad j \ne k 
\quad \text{and} \\ \EE \tau_j^2 \ = \ &O\left({1 \over n}\right) \quad \text{for} \quad j=1, \ldots, n.
\endaligned \tag7.1.1$$
For an unordered pair $1 \leq j \ne k \leq n$, let us define 
$$u_{\{j, k\}} =\root 3 \of{{1 \over 6} \zeta_{\{j, k\}} \left(1 -\zeta_{\{j, k\}} \right) \left(2 \zeta_{\{j, k\}} -1 \right) }
\left(\tau_j + \tau_k \right). $$
Then 
$$f=\sum_{\{j, k\}} u_{\{j, k\}}^3. $$ 
Similarly, let us define 
$$\aligned w_{\{j, k\}} = &\root 4 \of{{1 \over 24} \zeta_{\{j, k\}} \left(1-\zeta_{\{j, k\}}\right) \left| 6 \zeta_{\{j, k\}}^2 -6 \zeta_{\{j, k\}} +1\right|} \left(\tau_j +\tau_k \right) \quad \text{and} \\ 
\sigma_{\{j, k\}} =&\sgn \left(6 \zeta_{\{j, k\}}^2 - 6 \zeta_{\{j, k\}}+1 \right). \endaligned$$
Then
$$h=\sum_{\{j, k\}} \sigma_{\{j, k\}} w_{\{j, k\}}^4.$$
By (7.1.1) the random variables $u_{\{j, k\}}$ satisfy the conditions of Theorem 5.1 and hence 
the upper bound for $\mu = \EE f^2$ follows by Theorem 5.1. Similarly, by (7.1.1) the random 
variables $w_{\{j, k\}}$ satisfy the conditions of Theorem 6.1 and hence the upper bound for 
$|\nu| = \EE h$ and $\EE |h|$ follows by Part (1) of Theorem 6.1.
This concludes the proof of Part (2) of the theorem.

By Lemma 4.3, eigenvalues of $q$ are $\Omega(n)$ from which it follows that 
$$\left| \int_{{\Bbb R}^n \setminus {\Cal U}} e^{-q(t)} \ dt  \right| \ \leq \ \exp\bigl\{ -\Omega\left( \ln^2 n\right) \bigr\} \Xi. \tag7.1.2$$
From Theorem 5.1, we have 
$$\split &\left| \int_{{\Bbb R}^n} e^{-q(t)+if(t)} \ dt - \exp\left\{ -{\mu \over 2} \right\} \Xi \right| \ \leq \ 
\epsilon\ \Xi \\
&\qquad \qquad \text{provided} \quad  
n \geq \left( {1 \over \epsilon}\right)^{O(1)}, \endsplit$$
which, combined with (7.1.2), results in
$$\aligned &\left| \int_{\Cal U} e^{-q(t)+if(t)} \ dt - \exp\left\{ -{\mu \over 2} \right\} \Xi \right| \ \leq \ 
\epsilon\ \Xi \\
&\qquad \qquad \text{provided} \quad  
n \geq \left( {1 \over \epsilon}\right)^{O(1)}. \endaligned \tag7.1.3$$
By Part (2) of Theorem 6.1 and Chebyshev's inequality, we have 
$$\PP \Bigl\{ |h - \nu | > \epsilon\Bigr\} \ = \ O\left({1 \over \epsilon^2 n}\right), \tag7.1.4$$ 
while by Part (3) of Theorem 6.1, we have 
$$\PP\Bigl\{ h > \gamma(\delta)\Bigr\}  \ =O\left(\exp\left\{- n^{1/5}\right\}\right) \tag7.1.5$$
for some constant $\gamma(\delta)>0$.
In addition, 
$$|h(t)| =O\left(\ln^4 n\right) \quad \text{for} \quad t \in {\Cal U}. \tag7.1.6$$
Combining (7.1.4)--(7.1.6) and Part (2) of the theorem, we deduce from (7.1.2) and (7.1.3)
that 
$$\aligned
&\left| \int_{\Cal U} e^{-q(t) +h(t)} \ dt - \exp\left\{\nu\right\} \Xi \right| \ \leq \ \epsilon\  \Xi  
\\ &\qquad \qquad  \text{and} \\
&\left| \int_{\Cal U} e^{-q(t) +if(t)+ h(t)} \ dt - \exp\left\{ -{\mu \over 2} +\nu \right\} \Xi \right| 
\ \leq \ \epsilon \ \Xi \\
&\qquad \text{provided} \quad n \geq \left({1 \over \epsilon}\right)^{O(1)}. 
\endaligned \tag7.1.7$$

From the Taylor series expansion, cf. (3.2), we obtain
$$\split F(t)\ =\ &\exp\Bigl\{-q(t)+if(t)+h(t) +\rho(t) \Bigr\}, \quad \text{where} \\
|\rho(t)| \ =\ &O\left({\ln^5 n \over \sqrt{n}}\right) \quad \text{for} \quad 
t \in {\Cal U}. \endsplit $$
Therefore, for any $\epsilon >0$ we have 
$$\split &\Big||F(t)|-e^{-q(t)+h(t)}\Big|  \ \leq \ \epsilon\ e^{-q(t) + h(t)} 
\quad \text{and} \\ 
&\Big|F(t)-e^{-q(t)+if(t)+ h(t)}\Big| \leq \ \epsilon\ e^{-q(t) + h(t)}  \quad 
\text{for all} \quad t \in {\Cal U} \\
&\qquad \text{provided} \quad n \ \geq \ \left({1 \over \epsilon}\right)^{O(1)}.
\endsplit$$
The proof of Parts (3) and (4) now follows from (7.1.7) and Part (2).
 {\hfill \hfill \hfill} \qed
\enddemo

\head 8. Bounding the integral outside of the special points  \endhead 

We consider the integral representation of Corollary 3.3. Our goal is to show that the 
integral of $F(t)$ for $t$ outside of the neighborhood of the special points 
$$\left(0, \ldots, 0\right) \quad \text{and} \quad \left(\pm \pi, \ldots, \pm \pi \right)$$ 
is asymptotically negligible. 

In this section, we prove the following main result. 

\proclaim{(8.1) Theorem} Let us fix a number $0 < \delta \leq 1/2$ and let 
$D=\left(d_1, \ldots, d_n\right)$ be a $\delta$-tame degree sequence.
 Let us define subsets ${\Cal U}, {\Cal W} \subset \Pi$ by 
 $$\split {\Cal U}=&\Bigl\{\left(\tau_1, \ldots, \tau_n \right): \quad \left|\tau_j\right|\  \leq \ {\ln n \over \sqrt{n}} 
\quad \text{for} \quad j=1, \ldots, n\Bigr\}
 \quad \text{and} \\
{\Cal W}=&\Bigl\{ \left(\tau_1, \ldots, \tau_n\right): \quad \left| \tau_j - \sigma_j \pi \right| \ \leq \ 
{\ln n \over \sqrt{n}} \quad \text{for some} \quad \sigma_j =\pm 1\\ 
&\qquad \qquad \qquad \qquad \text{and all} \quad j=1, \ldots, n \Bigr\}. \endsplit$$
Then for any $\kappa>0$ 
$$\int_{\Pi \setminus ({\Cal U} \cup {\Cal W})} |F(t)| \ dt  \ \leq \ n^{-\kappa}  \int_{{\Cal U}} |F(t)| \ dt  $$
provided $n > \gamma(\delta, \kappa)$.
\endproclaim

The plan of the proof of Theorem 8.1 is as follows: first, using some combinatorial arguments we 
show that for any positive constant $\epsilon>0$ the integral is asymptotically negligible outside of the areas where $|\tau_j| \leq \epsilon$ for all $j$ or where $|\tau_j -\sigma_j \pi | \leq \epsilon$ for some $\sigma_j =\pm 1$ and all $j$. Then we note that $|F(t)|$ is log-concave for $t$ in a 
neighborhood of the origin and use a concentration inequality for log-concave measures.

We introduce the following metric $\rho$.

\subhead (8.2) Metric $\rho$ \endsubhead Let us define a function 
$\rho: {\Bbb R} \longrightarrow [0, \pi]$ as follows:
$$\rho(x)=\min_{k \in {\Bbb Z}} |x-2 \pi k|.$$
In words: $\rho(x)$ is the distance from $x$ to the nearest integer multiple of $2 \pi$. 
Clearly, 
$$\rho(-x)=\rho(x) \quad \text{and} \quad \rho(x+y) \leq \rho(x) +\rho(y)$$
for all $x, y \in {\Bbb R}$.

We will use that 
$$1-{1 \over 2} \rho^2(x) \ \leq \ \cos x \ \leq \ 1-{1 \over 5} \rho^2(x). \tag8.2.1$$

\subhead (8.3) The absolute value of $F(t)$ \endsubhead 
Let 
$$\alpha_{\{j, k\}} = 2 \zeta_{\{j, k\}} \left(1-\zeta_{\{j, k\}}\right) \quad \text{for all} \quad j \ne k.$$
If $D$ is $\delta$-tame, we have 
$$2 \delta^2 \ \leq \ \alpha_{\{j, k\}} \ \leq \ {1 \over 2} \quad \text{for all} \quad j \ne k. \tag8.3.1$$
We have 
$$|F(t)| =\left(\prod_{\{j, k\}} \left(1-\alpha_{\{j, k\}} + \alpha_{\{j, k\}} \cos \left(\tau_j + \tau_k \right) \right)
\right)^{1/2}. $$
For $1 \leq j \ne k \leq n$ let us define a function of $\tau \in {\Bbb R}$,
$$f_{\{j,k\}}( \tau)=\sqrt{1-\alpha_{\{j, k\}} +\alpha_{\{j, k\}} \cos \tau},$$
so 
$$|F(t)| =\prod_{\{j, k\}} f_{\{j, k\}}\left(\tau_j + \tau_k\right). \tag8.3.2$$
We note that 
$$f_{\{j, k\}}(0)=1.$$
It follows by (8.2.1) and (8.3.1) that for $\epsilon>0$
$$\aligned &f_{\{j, k\}}(x) \ \leq \ \exp\bigl\{ -\gamma(\delta, \epsilon)\bigr\} f_{\{j, k\}}(y)\\
&\qquad \text{provided} \quad \rho(x) \geq 2 \epsilon \quad \text{and} \quad \rho(y) \leq \epsilon
\endaligned \tag8.3.3$$
for some $\gamma(\epsilon, \delta)>0$.
Furthermore,
$${d^2 \over d \tau^2} \ln f_{\{j, k\}}(\tau)=-{\alpha_{\{j, k\}}\left(\alpha_{\{j, k\}} + \cos \tau - \alpha_{\{j, k\}} \cos \tau\right)  \over 2\left(1-\alpha_{\{j, k\}} +\alpha_{\{j, k\}} \cos \tau\right)^2}. $$
In particular, by (8.3.1)
$${d^2 \over d \tau^2} \ln f_{\{j, k\}}(\tau) \ \leq \ -{\delta^2 \over 2}
\quad \text{for} \quad -{\pi \over 3} \ \leq \ \tau 
\ \leq \ {\pi \over 3}$$
and hence $\ln f_{\{j, k\}}$ is strictly concave on the interval $[-\pi/3, \pi/3]$:
$$\aligned &\ln f_{\{j, k\}}(x)  +\ln f_{\{j, k\}}(y) - 2\ln f_{\{j, k\}}\left({x+y \over 2}\right) \ \leq \ -{\delta^2 \over 8}|x-y|^2 \\
&\qquad \text{for all} \quad x, y \in [-\pi/3, \pi/3].\endaligned \tag8.3.4$$
\bigskip
In what follows, we fix a particular parameter $\epsilon >0$. All implied constants in the
``$O$'' and ``$\Omega$'' notation below may depend only on the parameters $\delta$ and $\epsilon$.
We say that $n$ is {\it sufficiently large} if $n \geq \gamma(\delta, \epsilon)$ for some 
constant $\gamma(\delta, \epsilon)>0$.
\bigskip
Our first goal is to show that only the points $t \in \Pi$ for which the inequality 
$\rho\left(\tau_j +\tau_k \right) \leq \epsilon$ holds for an overwhelming majority of pairs $\{j, k\}$ contribute
significantly to the integral of $|F(t)|$ on $\Pi$.
\proclaim{(8.4) Lemma} For $t \in \Pi$, $t=\left(\tau_1, \ldots, \tau_n \right)$, and $\epsilon>0$
let us define a set 
$K(t, \epsilon) \subset \{1, \ldots, n\}$ consisting of the indices $k$ such that 
$$\rho\left( \tau_j +\tau_k \right)\ \leq \ \epsilon$$
for more than $n/2$ distinct indices $j$.
Let $\overline{K(t, \epsilon)}=\{1, \ldots, n\} \setminus K(t, \epsilon)$.
Then
\roster
\item
$$|F(t)| \ \leq\ \exp\Bigl\{-\gamma(\delta) \epsilon^2 n |\overline{K(t, \epsilon)}|  \Bigr\}
\quad \text{for some} \quad \gamma(\delta)>0;$$ 
\item
$$\rho\left(\tau_{k_1} -\tau_{k_2} \right) \ \leq \ 2\epsilon  
\quad \text{for all} \quad k_1, k_2 \in K(t, \epsilon);$$
\item Suppose that $|K(t, \epsilon)| > n/2$. Then 
$$\rho\left(\tau_{k_1} +\tau_{k_2} \right) \ \leq \ 3\epsilon \quad \text{for all} \quad k_1, k_2 \in 
K(t, \epsilon).$$
\endroster
\endproclaim
\demo{Proof} For every $k \in \overline{K(t, \epsilon)}$ there are at least $(n-2)/2$ distinct $j \ne k$ for which 
$$\rho\left(\tau_j + \tau_k\right) \ > \ \epsilon\tag8.4.1$$ and so by (8.2.1) we have 
$$\cos\left(\tau_j +\tau_k \right) \ \leq \ 1- {1 \over 5} \epsilon^2. $$ Since there are at least 
$|\overline{K(t, \epsilon)}|(n-2)/4$ pairs $\{j, k\}$ for which (8.4.1) holds,
the proof of Part (1) follows from (8.3.2) and (8.3.3).

For any $k_1, k_2 \in K$ there is $j \in \{1, \ldots, n\}$ such that 
$$\rho\left(\tau_j + \tau_{k_1} \right), \ \rho\left(\tau_j + \tau_{k_2} \right) \ \leq \ \epsilon.$$
Therefore,
$$\split \rho\left( \tau_{k_1}-\tau_{k_2} \right) =&\rho\left(\tau_j +\tau_{k_1}  -\tau_{k_2}-\tau_j\right) 
\leq \rho\left(\tau_j +\tau_{k_1}\right) + \rho\left(-\tau_j-\tau_{k_2}  \right) \\=
&\rho\left(\tau_j+ \tau_{k_1} \right) +\rho\left(\tau_j + \tau_{k_2} \right) \ \leq \ 2\epsilon
\endsplit$$
and Part (2) follows.

Let us choose a $k_1 \in K(t, \epsilon)$. Since $|K(t, \epsilon)| > n/2$, there is a $j \in K(t, \epsilon)$
such that 
$$\rho\left(\tau_j + \tau_{k_1} \right) \ \leq \ \epsilon.$$
By Part (2), for any $k_2 \in K(t, \epsilon)$ we have 
$$\split \rho\left(\tau_{k_1} + \tau_{k_2} \right) \ =\ &\rho\left(-\tau_j+\tau_{k_1} +\tau_j + \tau_{k_2}\right) 
\ \leq \ \rho\left(\tau_j +\tau_{k_1} \right) + \rho\left(-\tau_j +\tau_{k_2} \right) \\
\leq \ &3 \epsilon  \endsplit$$
and Part (3) follows.
{\hfill \hfill \hfill} \qed
\enddemo 

\proclaim{(8.5) Corollary} For an $\epsilon >0$  let us define 
a set $V(\epsilon) \subset \Pi$ consisting of the points $t \in \Pi$ such that 
$$\overline{K(t, \epsilon)} \ \leq \  \ln^2 n, $$
where $\overline{K(t, \epsilon)}$ is defined in Lemma 8.4. Then
$$\int_{\Pi \setminus V(\epsilon)} |F(t)| \ dt \ \leq \ n^{-n} \int_{\Pi} |F(t)| \ dt$$
provided $n \geq \gamma(\delta, \epsilon)$ for some constant $\gamma(\delta, \epsilon)>0$.
\endproclaim
\demo{Proof} By Parts (1)--(3) of Theorem 7.1
we have 
$$\int_{\Pi} |F(t)| \ dt \ \geq \ \Omega\left(n^{-n/2}\right).$$
The proof now follows by Part (1) of Lemma 8.4.
{\hfill \hfill \hfill} \qed 
\enddemo

Next, we show that only the points $t \in \Pi$ such that $\rho\left(\tau_j +\tau_k \right) \leq \epsilon$
for {\it all} $1 \leq j, k \leq n$ contribute substantially to the integral of $|F(t)|$ on $\Pi$.

\proclaim{(8.6) Lemma} For $\epsilon>0$ let us define a set $X(\epsilon) \subset \Pi$,
$$X(\epsilon)=\Bigl\{t \in \Pi, \ t=\left(\tau_1, \ldots, \tau_n\right): \ 
\rho\left(\tau_j +\tau_k \right) \leq \epsilon \quad \text{for all} \quad j, k\Bigr\}.$$
Then
$$\int_{\Pi \setminus X(\epsilon)} |F(t)| \ dt \ \leq \ \exp\bigl\{-\gamma_1(\delta, \epsilon) n \bigr\}
\int_{\Pi} |F(t)| \ dt$$
for all $n \geq \gamma_2(\delta, \epsilon)$ for some constants 
$\gamma_1(\delta, \epsilon), \gamma_2(\delta, \epsilon)>0$.
\endproclaim
\demo{Proof} Let us consider the set $V(\epsilon/60) \subset \Pi$ and $n$ large enough so that the conclusion 
of Corollary 8.5 holds, that is, the integral of $|F(t)|$ over $\Pi \setminus V(\epsilon/60)$ is 
asymptotically negligible. For a set $A \subset \{1, \ldots, n\}$ such that 
$$|\overline{A}| \ \leq \  \ln^2 n,$$
let us define a set $P_A \subset \Pi$ (we call it a {\it piece}) such that 
$$\rho\left(\tau_j -\tau_k \right) \ \leq \ \epsilon/30 \quad \text{and} \quad \rho\left(\tau_j +\tau_k \right) 
\ \leq \ \epsilon/20 \quad \text{for all} \quad j, k \in A.$$
If $n$ is large enough,  by Lemma 8.4
for every $t \in V(\epsilon/60)$ we can choose $A=K(t, \epsilon)$, so we have 
$$V(\epsilon/60) \ \subset \ \bigcup_{A: \ |\overline{A}|\leq \ln^2 n} P_A. \tag8.6.1$$
Our next goal is to show that the integral of $|F(t)|$ over $P_A \setminus X(\epsilon)$ is 
negligible compared to the integral of $|F(t)|$ over $P_A$.

Let us choose a point $t \in P_A \setminus X(\epsilon)$. Thus we have 
$$\rho\left(\tau_{i_0} + \tau_{j_0} \right) > \epsilon \quad \text{for some} \quad i_0, j_0.$$
Let us choose any $k_0 \in A$. Then
$$\rho\left(\tau_{i_0} +\tau_{j_0}\right) =\rho\left(\tau_{i_0} + \tau_{j_0} + \tau_{k_0} -\tau_{k_0}\right)
\leq \rho\left(\tau_{i_0} +\tau_{k_0}\right) +\rho\left(\tau_{j_0} -\tau_{k_0}\right).$$
Hence we have either 
$$\rho\left(\tau_{i_0} + \tau_{k_0} \right) > \epsilon/2 \quad \text{or} \quad 
\rho\left(\tau_{j_0}-\tau_{k_0}\right) > \epsilon/2.$$
In the first case, for every $k \in A$ we have 
$$\split \rho\left(\tau_{i_0} +\tau_{k_0} \right)=&\rho\left(\tau_{i_0} +\tau_{k_0} +\tau_k -\tau_k \right)  \ \leq \ \rho\left(\tau_{i_0} +\tau_k \right) + \rho\left(\tau_{k_0} -\tau_k \right) \\
&\leq \rho\left(\tau_{i_0} +\tau_k \right) +\epsilon/30, \endsplit$$
from which 
$$\rho\left(\tau_{i_0} +\tau_k \right) \ \geq \ \epsilon/2 -\epsilon/30 = 7\epsilon/15.$$
In the second case, for every $k \in A$, we have 
$$\split\rho\left(\tau_{j_0} - \tau_{k_0} \right) =
&\rho\left(\tau_{j_0} -\tau_{k_0} + \tau_k - \tau_k \right) \ \leq \ 
\rho\left(\tau_{j_0} +\tau_k\right) + \rho\left(-\tau_{k_0} -\tau_k\right)\\ =
&\rho\left(\tau_{j_0} +\tau_k \right) + \rho\left(\tau_{k_0} +\tau_k \right) \ \leq \ 
\rho\left(\tau_{j_0} + \tau_k \right) + \epsilon/20, \endsplit$$
from which 
$$\rho\left(\tau_{j_0} +\tau_k \right) > \epsilon/2 - \epsilon/20 = 9\epsilon/20.$$
In either case, for any $t \in P_A\setminus X(\epsilon)$, $t=\left(\tau_1, \ldots, \tau_n\right)$, there exists an index $i \notin A$ such that
$$\rho\left(\tau_i + \tau_k \right) > 0.45 \epsilon \quad \text{for all} \quad k \in A.$$
For $i \notin A$, we define 
$$Q_{A, i}=\Bigl\{t \in P_A: \ \rho\left(\tau_i +\tau_k \right) > 0.45 \epsilon \quad \text{for all} \quad k \in 
A \Bigr\}. \tag8.6.2$$ 
Hence 
$$P_A \setminus X(\epsilon)  \ \subset \ \bigcup_{i \in \overline{A}} Q_{A, i}. $$
Given a point $t \in P_A$, we obtain another point in $P_A$ if we arbitrarily change 
the coordinate $\tau_i \in [-\pi, \pi]$ for $i \in \overline{A}$. We obtain a {\it fiber} $E \subset P_A$ 
if we fix all other coordinates and let $\tau_i \in [-\pi, \pi]$ vary. Geometrically, each fiber $E$ 
is an interval of length $2\pi$. Let us construct a set $I \subset E$ as follows. We choose 
an arbitrary $k_0 \in A$ and let $\tau_i$ vary in such a way that $\rho\left(\tau_{k_0} +\tau_i \right) 
\leq 0.05 \epsilon$. Geometrically, $I$ is an interval of length $0.1 \epsilon$ or a union of 
two non-overlapping intervals of total length $0.1 \epsilon$. Moreover, 
$$\rho\left(\tau_k +\tau_i \right) \ \leq \ \rho\left(\tau_{k_0} + \tau_i \right) + \rho\left( \tau_k -\tau_{k_0}\right) \ \leq \ 0.05 \epsilon + 0.05 \epsilon =0.1 \epsilon$$ 
for all $k \in A$ and all $\tau \in I$. 

Using (8.3.2) and (8.3.3), we conclude from (8.6.2) that for any $t \in Q_{A, i} \cap E$ and for any $s \in P_A \cap I$ we have 
$$|F(t)| \ \leq \exp\bigl\{-\Omega(n) \bigr\} |F(s)|$$ 
provided $n$ is large enough.
Therefore,
$$\int_{E \cap Q_{A, i}} |F(t)| \ dt \ \leq \ \exp\bigl\{-\Omega(n)\bigr\} \int_E |F(t)| \ dt$$
for all sufficiently large $n$.

Integrating over all fibers $E$, we establish that 
$$\int_{Q_{A, i}} |F(t)| \ dt \ \leq \ \exp\bigl\{-\Omega(n)\bigr\} \int_{P_A} |F(t)| \ dt$$
for all sufficiently large $n$. 
Since the number of different subsets $A \subset \{1, \ldots, n\}$ with $|\overline{A}| \leq  \ln^2 n$ in 
(8.6.1)
does not exceed $\exp\bigl\{O\left( \ln^3 n\right) \bigr\}$, the proof follows.
{\hfill \hfill \hfill}\qed
\enddemo

Next, we prove that only the points in the neighborhood of the origin or the corners of $\Pi$ 
contribute significantly to the integral of $|F(t)|$ over $\Pi$.

\proclaim{(8.7) Lemma} For $0 < \epsilon < 1$, let $X(\epsilon)$ be the set defined in Lemma 8.6.
Let us define $Y(\epsilon), Z(\epsilon) \subset \Pi$ by 
$$\split Y(\epsilon)=&\Bigl\{t \in \Pi: \ t=\left(\tau_1, \ldots, \tau_n \right), \quad 
|\tau_i| \leq \epsilon/2 \quad \text{for} \quad i=1, \ldots, n \Bigr\} \quad \text{and} \\
Z(\epsilon)=&\Bigl\{t \in \Pi: \ t=\left(\tau_1, \ldots, \tau_n \right), \quad 
|\tau_i -\sigma_i \pi| \leq \epsilon/2 \quad \text{for some} \quad \sigma_i=\pm 1\\
& \qquad \quad \text{and all} \quad i=1, \ldots, n \Bigr\}. \endsplit$$
Then
$$X(\epsilon) \ = \ Y(\epsilon) \cup Z(\epsilon) \quad \text{and} \quad
 Y(\epsilon) \cap Z(\epsilon)=\emptyset.$$
Moreover, 
$$\int_{Y(\epsilon)} |F(t)| \ dt =\int_{Z(\epsilon)} |F(t)| \ dt.$$
\endproclaim
\demo{Proof} Let us pick a point $t \in X(\epsilon)$. Then for each $k$ we have 
$\rho\left(2 \tau_k \right) \leq \epsilon$ and hence either
$$|\tau_k| \leq \epsilon/2 \quad \text{or} \quad |\tau_k -\pi| \leq \epsilon/2 \quad \text{or} \quad
|\tau_k +\pi| \leq \epsilon/2.$$
Since $\rho\left(\tau_k +\tau_j \right) \leq \epsilon$ for all $k, j$, we conclude that if 
$|\tau_k| \leq \epsilon/2$ for some $k$ then $|\tau_k| \leq \epsilon/2$ for all $k$. 
Hence $X(\epsilon) \ \subset  \ \bigl(Y(\epsilon) \cup Z(\epsilon)\bigr)$. The inclusion 
$\bigl(Z(\epsilon)  \cup Y(\epsilon)\bigr) \ \subset \ X(\epsilon)$ is obvious.  Since $\epsilon<1$, we have 
$Y(\epsilon) \cap Z(\epsilon) =\emptyset$.

The set $Z(\epsilon)$ is a union of $2^n$ pairwise disjoint {\it corners}, where each corner 
is determined by a choice of the interval $[-\pi, -\pi +\epsilon/2]$ or $[\pi-\epsilon/2, \pi]$
for each coordinate $\tau_i$. The transformation 
$$\tau_k \longmapsto \cases \tau_k +\pi &\text{if\ } \tau_k \in [-\pi, -\pi+\epsilon/2] \\
\tau_k -\pi &\text{if \ } \tau_k \in [\pi-\epsilon/2, \pi] \endcases$$
is a volume-preserving transformation which maps $Z(\epsilon)$ onto $X(\epsilon)$ and does not 
change the value of $|F(t)|$.
{\hfill \hfill \hfill} \qed
\enddemo

Finally, we will use that $|F(t)|$ is strictly log-concave on the set $Y(1)$. For Euclidean space $V$ with the norm 
$\| \cdot \|$, a point $x \in V$ and a closed set $A \subset V$ we define the distance
$$\dist(x, A)=\min_{y \in A} \|x-y\|.$$
We will need the 
following concentration inequality for strictly log-concave measures.
\proclaim{(8.8) Theorem} Let $V$ be Euclidean space with the norm $\| \cdot \|$, let 
$B \subset V$ be a convex body and let us consider a probability measure supported on 
$B$ with density $e^{-u}$, where $u: B \longrightarrow {\Bbb R}$ is a function satisfying 
$$u(x) +u(y) -2u\left({x+y \over 2} \right) \geq c \|x-y\|^2 \quad \text{for all} \quad x, y \in B$$
and some constant $c >0$. 

Let $A \subset B$ be a closed subset such that $\PP(A) \geq 1/2$. Then, for $r \geq 0$, we have 
$$\PP\Bigl\{x \in B: \quad \dist(x, A) \geq r \Bigr\} \leq 2 e^{-cr^2}.$$
\endproclaim
\demo{Proof} See Section 2.2 of
\cite{Le01} and Theorem 8.1 and its proof in \cite{Ba97}.
{\hfill \hfill \hfill} \qed
\enddemo

\proclaim{(8.9) Lemma} Let $Y(1) \subset \Pi$ be the set defined by
$$Y(1)=\Bigl\{t \in \Pi, \ t=\left( \tau_1, \ldots, \tau_n \right): \quad |\tau_i| \ \leq \ {1 \over 2} \quad 
\text{for} \quad i=1, \ldots, n\Bigr\}.$$
Let ${\Cal U} \subset \Pi$ be the set 
$${\Cal U}=\Bigl\{t \in \Pi,\ t=\left(\tau_1, \ldots, \tau_n \right): \quad |\tau_i| \ \leq \ {\ln n \over \sqrt{n}} 
\quad \text{for} \quad i=1, \ldots, n \Bigr\}.$$
Then for any $\kappa >0$ we have 
$$\int_{Y(1) \setminus {\Cal U}}|F(t)| \ dt \ \leq \ n^{-\kappa} \int_{Y(1)} |F(t)| \ dt,$$
provided $n \geq \gamma(\delta, \kappa)$ is big enough.
{\hfill \hfill \hfill} \qed
\endproclaim
\demo{Proof} Let us consider the probability measure on $Y(1)$ with density proportional 
to $|F(t)|$. Let us consider a map $M: {\Bbb R}^n \longrightarrow {\Bbb R}^{n \choose 2}$, where the coordinates 
of ${\Bbb R}^{n \choose 2}$ are indexed by unordered pairs $\{j, k\}$ and 
$$M_{\{j, k\}} \left(\tau_1, \ldots, \tau_n \right) =\tau_j +\tau_k.$$
By Lemma 4.3, 
$$\| M (t)\|^2 \ \geq \ (n-2) \|t\|^2 \quad \text{for all} \quad t \in {\Bbb R}^n.$$
Since from Section 8.3,
$$\ln |F(t)| ={1 \over 2} \sum_{\{j, k\}} \ln f_{\{j, k\}}\left(\tau_j +\tau_k \right),$$
it follows by (8.3.4) that for any constant $a$ 
and 
$$u(t)=-\ln |F(t)| + a,$$ 
we have 
$$u\left(t_1\right) + u\left(t_2 \right) - 2 u\left( {t_1+t_2 \over 2} \right) \ \geq \ \gamma(\delta)n 
\|t_1 -t_2\|^2 \quad \text{for all} \quad t_1, t_2 \in Y(1)$$
and some constant $\gamma(\delta)>0$.
We choose $a$ so that $e^{-u}$ is a probability density on $Y(1)$. 

We apply Theorem 8.8 
with $c=\gamma(\delta) n$. For $k=1, \ldots, n$, let $A_k^- \subset Y(1)$ be the set of points 
with $\tau_k \leq 0$ and let $A_k^+ \subset Y(1)$ be the set of points with $\tau_k \geq 0$. 
Since both $Y(1)$ and the probability measure are invariant under the symmetry 
$t \longmapsto -t$, we have 
$$\PP\left(A_k^-\right)=\PP\left(A_k^+\right)={1 \over 2} \quad \text{for} \quad k=1, \ldots, n.$$
Therefore, by Theorem 8.8, all but a $n^{-\kappa}$ fraction of all points in $Y(1)$ lie within 
a distance of $\ln n/\sqrt{n}$ from each of the sets $A_k^-$ and $A_k^+$, provided 
$n$ is large enough.
{\hfill \hfill \hfill} \qed
\enddemo
\subhead (8.10) Proof of Theorem 8.1 \endsubhead 
For $0 < \epsilon < 1$ let us define the set $X(\epsilon)$ as in Lemma 8.6 and the sets 
$Y(\epsilon)$ and $Z(\epsilon)$ as in Lemma 8.7. In particular,
$${\Cal U}=Y\left({2 \ln n \over \sqrt{n}} \right) \quad \text{and} \quad {\Cal W} = Z\left({2 \ln n \over \sqrt{n}} \right).$$

By Lemma 8.6, for any $\kappa>0$ we have 
$$\int_{\Pi \setminus X(1)} |F(t)| \ dt \ \leq \ n^{-\kappa} \int_{\Pi} |F(t)| \ dt$$
for all sufficiently large $n$, so that the integral outside of $X(1)$ is asymptotically negligible.

By Lemma 8.7, $X(1)=Y(1) \cup Z(1)$ with $Y(1) \cap Z(1)=\emptyset$ and 
$$\int_{Y(1)}  |F(t)| \ dt = \int_{Z(1)} |F(t)| \ dt \quad \text{and} \quad 
\int_{\Cal U} |F(t)| \ dt =\int_{\Cal W} |F(t)| \ dt. \tag8.10.1$$
By Lemma 8.9,
$$\int_{Y(1) \setminus {\Cal U}} |F(t)| \ dt \ \leq \ n^{-\kappa} \int_{Y(1)} |F(t)| \ dt$$
for all sufficiently large $n$, so that the integral over $Y(1) \setminus {\Cal U}$ is asymptotically negligible. 
By (8.10.1), the integral over $Z(1) \setminus {\Cal W}$  is asymptotically negligible. 
The proof now follows.
{\hfill \hfill \hfill} \qed

\head 9. Proof of Theorem 1.4 \endhead 

By Corollary 3.3, we have the integral representation for the number $|G(D)|$ of graphs:
$$|G(D)| ={e^{H(z)} \over (2 \pi)^n} \int_{\Pi} F(t) \ dt.$$
Let us define subsets ${\Cal U}, {\Cal W} \subset \Pi$ as in Theorem 8.1. Let us consider the transformation 
${\Cal W} \longrightarrow {\Cal U}$, 
$$\tau_k \longmapsto \cases \tau_k + \pi &\text{if\quad } 
-\pi  \ \leq \ \tau_k \ \leq \ -\pi + {\ln n \over \sqrt{n}} \\ 
\tau_k - \pi &\text{if\quad } \pi - {\ln n \over \sqrt{n}} \ \leq \ \tau_k \ \leq \pi \endcases 
\qquad \quad \text{for} \quad k=1, \ldots, n.$$
As in the proof of Lemma 8.7, this is a measure-preserving transformation which maps 
${\Cal W}$ onto ${\Cal U}$. Since $d_1 + \ldots  +d_n$ is even, the transformation does not change the value of 
$F(t)$ (if $d_1 + \ldots  + d_n$ is odd, the transformation changes the sign of $F(t)$).
Hence 
$$\int_{\Cal U} F(t) \ dt = \int_{\Cal W} F(t) \ dt. \tag9.1$$
By Theorem 7.1, the integrals of $F(t)$ and $|F(t)|$ over ${\Cal U}$ have the same order of magnitude, that is, 
$$\int_{\Cal U} |F(t)| \ dt \ \leq \ \gamma(\delta) \left| \int_{\Cal U} F(t) \ dt \right| $$
for some constant $\gamma(\delta) >1$. Therefore, from Theorem 8.1, the integral outside of 
${\Cal U} \cup {\Cal W}$ is asymptotically negligible, so that for any $\kappa >0$ and all sufficiently large 
$n \ \geq \gamma(\delta, \kappa)$, we have 
$$\int_{\Pi \setminus ({\Cal U} \cup {\Cal W})} |F(t)| \ dt \ \leq \ n^{-\kappa} \left|  \int_{\Cal U} F(t) \ dt \right|.$$
The proof now follows by Parts (2) and (4) of Theorem 7.1, identity (9.1) and the formula 
$$\Xi \ =\ \int_{{\Bbb R}^n} e^{-q(t)} \ dt \ = \ {(2 \pi)^{n/2} \over \sqrt{\det Q}}.$$
{\hfill \hfill \hfill} \qed

\head 10. Proof of Theorem 1.6   \endhead 

The proof is very similar to that of Theorem 3 of \cite{Ba10},  which 
deals with a similar situation in the case of bipartite graphs. All implicit constant in 
the ``$O$'' and ``$\Omega$''-notation below may depend only on the parameter $\delta>0$.

For pairs $1 \leq j \ne k \leq n$, let $x_{\{j, k\}}$ be independent Bernoulli random variables such that 
$$\PP\bigl\{x_{\{j, k\}} =1 \bigr\} =\zeta_{\{j, k\}} \quad \text{and} \quad 
\PP\bigl\{x_{\{j, k\}} =0 \bigr\} =1-\zeta_{\{j, k\}}.$$
As is implied by Theorem  5 of \cite{BH10}, the probability mass function of the random vector 
$X=\left(x_{\{j, k\}}\right)$
is constant on the integer points of ${\Cal P}(D)$ and is equal to $e^{-H(z)}$ at each 
$G\in G(D)$.

Let us define
$$\sigma_S(X)=\sum_{\{j, k\} \in S} x_{\{j, k\}}.$$
Then 
$$\split &\PP \Bigl\{G \in G(D): \quad \sigma_S(G) \ \leq \ (1-\epsilon) \sigma_S(z) \Bigr\} 
\\=&\qquad{\PP \Bigl\{X : \quad \sigma_S(X) \ \leq \ (1-\epsilon) \sigma_S(z) \quad \text{and} \quad 
X \in G(D) \Bigr\} \over \PP\Bigl\{X: \quad X \in G(D) \Bigr\}} \endsplit$$ and,
similarly,
$$\split &\PP \Bigl\{G \in G(D): \quad \sigma_S(G) \ \geq \ (1+\epsilon) \sigma_S(z) \Bigr\} 
\\=&\qquad{\PP \Bigl\{X : \quad \sigma_S(X) \ \geq \ (1+\epsilon) \sigma_S(z) \quad \text{and} \quad 
X \in G(D) \Bigr\} \over \PP\Bigl\{X: \quad X \in G(D) \Bigr\}} \endsplit$$
Applying Theorem 1.4 and Parts (1) and (2) of Theorem 7.1, we get
$$\PP \Bigl\{X: \quad X \in G(D) \Bigr\} =e^{-H(z)} |G(D)|  \ \geq \ n^{-O(n)}.$$
On the other hand, standard large deviation inequalities
for sums of bounded independent random variables (see, for example, Corollary 5.3 of \cite{McD89})
imply that 
$$\PP \Bigl\{X : \quad \sigma_S(X) \ \geq \ (1+\epsilon) \sigma_S(z)\Bigr\} \ \leq \  \exp\left\{ -\Omega\left( n  \ln^2 n \right)\right\}$$
and, similarly, 
$$\PP \Bigl\{X : \quad \sigma_S(X) \ \leq \ (1-\epsilon) \sigma_S(z)\Bigr\} \ \leq \ 
\exp\left\{ -\Omega\left( n  \ln^2 n\right) \right\}$$
 and the proof follows.
 {\hfill \hfill \hfill} \qed

\head 11. Counting bipartite graphs \endhead

Here we list some modifications needed to establish the asymptotic formula (2.5.4). 
We adhere to the notation of Section 2.5.

As in  Corollary 3.3, we represent the number 0-1 matrices with row sums 
$R=\left(r_1, \ldots, r_m \right)$ and column sums $C=\left(c_1, \ldots, c_n \right)$ as 
an integral. Let us define 
$$\split F(s, t)=&\exp\left\{-i\sum_{j=1}^m r_j \sigma_j -i \sum_{k=1}^n c_k \tau_k \right\} 
\prod \Sb 1 \leq j \leq m \\ 1 \leq k \leq n \endSb \left(1-\zeta_{jk} + \zeta_{jk} e^{i(\sigma_j + \tau_k)} \right)
\\ & \qquad \qquad  \text{for} \quad (s, t)=\left(\sigma_1, \ldots, \sigma_m; \tau_1, \ldots, \tau_n \right)
\endsplit$$
and let $\Pi \subset {\Bbb R}^{m+n}$ be the parallelepiped 
$$\split \Pi=&\Bigl\{ \left(\sigma_1, \ldots, \sigma_m; \tau_1, \ldots \tau_n \right): \quad 
-\pi \ \leq \ \sigma_j, \tau_k \leq \pi \\ &\qquad \qquad  \text{for} \quad j=1, \ldots, m;\ k=1, \ldots, n \Bigr\}.
\endsplit$$
Let $\Pi_0 \subset \Pi$ be the facet of $\Pi$ defined by the equation $\tau_n=0$.

Since the constraints are not independent (the sum of all row sums is equal to the 
sum of all column sums), we can drop one of the constraints and represent the desired number $|R, C|$ as an integral over 
$\Pi_0$,
$$|R, C| = {e^{H(z)} \over (2\pi)^{m+n-1}} \int_{\Pi_0} F(s, t) \ ds\, dt,$$
cf. Section 2 of \cite{BH09a}.

Let ${\Cal U} \subset \Pi$ be the neighborhood of the origin,
$$\split {\Cal U} =&\Bigl\{ \left(\sigma_1, \ldots, \sigma_m; \tau_1, \ldots, \tau_n \right): \quad 
\left| \sigma_j \right|, \ \left| \tau_k \right| \ \leq {\ln n \over \sqrt{n}} \\
&\qquad \qquad \text{for} \quad j=1, \ldots, m;\ k=1, \ldots, n \Bigr\}\endsplit$$ 
and let ${\Cal U}_0$ be the intersection of ${\Cal U}$ with the hyperplane $\tau_n=0$. 
We prove that the integral 
$$\int_{\Pi_0 \setminus {\Cal U}_0} |F(s, t)| \ ds\, dt$$
is asymptotically negligible relative to the integral 
$$\int_{{\Cal U}_0} |F(s, t)| \ ds\, dt. \tag11.1$$ 
The proof is a modification of that of Theorem 8.1 (note that here we don't have another 
set ${\Cal W} \subset \Pi$ contributing large values of $|F(s,t)|$) and very similar to that of 
Theorem 7.1 of \cite{BH09a}. A different line of proof can be inferred from \cite{BH09b}.

Our next goal is to evaluate asymptotically as $m, n \longrightarrow +\infty$ the integral 
$$\int_{{\Cal U}_0} F(s, t)\ ds\, dt.\tag11.2$$
In particular, we need to show that (11.2) and (11.1) have about the same order of magnitude.
From (3.2), we can write the expansion 
$$F(s, t) = \exp\bigl\{-q(s, t) + i f(s, t) +h(s, t) \bigr\}\Bigl(1+o(1)\Bigr) \quad \text{for} \quad (s, t) \in {\Cal U}, $$
as $m, n \longrightarrow +\infty$, 
where $q, f$ and $h$ are as defined by formulas (2.5.1) and (2.5.3) respectively. 

We note that $q(s, t)$ is {\it not} strictly positive definite, since its kernel is spanned by the vector 
$u \in {\Bbb R}^{m+n}$ defined by (2.5.2). However, the restriction of $q$ onto any hyperplane 
$L \subset {\Bbb R}^{m+n}$ which does not contain $u$ is strictly positive definite and allows us 
to define the Gaussian probability measure in $L$ with density proportional to $e^{-q}$. 
It is easy to prove (see Lemma 3.1 of \cite{BH09a}) that the expectation of any polynomial in 
the sums $\sigma_j + \tau_k$ does not depend on the choice of $L$. 
To evaluate (11.2), we need to show that asymptotically
$$\EE \exp\bigl\{ i f + h \bigr\} = \exp\left\{ -{1 \over 2} \EE f^2 +  \EE h  \right\} \Bigl(1 +o(1)\Bigr),$$
if we choose the hyperplane $L$ defined by the equation $\tau_n=0$. However, since the 
expectation on the left hand side does not depend on the choice of the hyperplane $L$, we can choose 
$L$ in such a way that 
$$\split \left| \EE \tau_j \tau_k \right|, \ \left| \EE \sigma_j \sigma_k \right| \ = \ &O \left({1 \over mn} \right)
\quad \text{provided} \quad j \ne k \\
\left| \EE \sigma_j \tau_k \right| \ = \ &O\left({1 \over mn}\right) \quad \text{for all} \quad j, k \quad \text{and}
\\ 
\EE \sigma_j^2, \ \EE \tau_k^2 \ = \ &O\left({1 \over m+n} \right) \quad \text{for all} \quad j, k. 
\endsplit \tag11.3$$
As is shown in \cite{BH09a} (see Theorem 3.2 there), to ensure (11.3), one has to choose 
$L$ defined by the equation 
$$\split &\sum_{j=1}^m \alpha_j \sigma_j =\sum_{k=1}^n \beta_k \tau_k, \quad \text{where} \\
&\quad \alpha_j =\sum_{k=1}^n \left(\zeta_{jk}-\zeta_{jk}^2\right) \quad \text{and} \quad
\beta_k=\sum_{j=1}^m \left(\zeta_{jk}-\zeta_{jk}^2 \right). \endsplit$$
The proof then proceeds as in Theorem 1.4.

\head 12. Proof of Theorem 2.1 \endhead 

In what follows, it is convenient to define the polytope ${\Cal P}(D) \subset {\Bbb R}^{{n \choose 2}}$
for positive, not necessarily integer, sequences $D=\left(d_1, \ldots, d_n \right)$. 
Recall that ${\Cal P}(D)$ consists of the vectors $\left(\xi_{\{j, k\}} \right)$ for 
$1 \leq j \ne k \leq n$ such that 
$$\sum_{j: \ j \ne k} \xi_{\{j, k\}} = d_k \quad \text{for} \quad k=1, \ldots, n$$ 
and 
$$0 \ \leq \xi_{\{j, k\}} \ \leq 1 \quad \text{for} \quad 1 \leq j \ne k \leq n.$$
We say that ${\Cal P}(D)$ has a {\it non-empty interior} if there is a point 
$y \in {\Cal P}(D)$, $y=\left(\eta_{\{j, k\}}\right)$, such that 
$$0 \ < \ \eta_{\{j, k\}} \ < \ 1 \quad \text{for all} \quad 1 \leq j\ne k \leq n.$$

The following two lemmas are probably known in greater generality, but since we are 
unable to provide a precise reference, we prove only the parts we need to obtain Theorem 2.1.

 \proclaim{(12.1) Lemma} Let $D=\left(d_1, \ldots, d_n \right)$ be a sequence of positive rational numbers such that 
$$d_1 \geq \ldots \geq d_n$$
and the Erd\H{o}s-Gallai conditions
 $$\sum_{i=1}^k d_i \ \leq \ k(k-1)+  \sum_{i=k+1}^n \min \left\{k, d_i \right\} 
\quad \text{for} \quad k=1, \ldots, n$$
are satisfied. Then the polytope ${\Cal P}(D)$ is non-empty.
\endproclaim
\demo{Proof}  Let $q$ be a positive integer
such that $qd_i$ are even integer for $i=1, \ldots, n$. Clearly, ${\Cal P}(D)$ is non-empty if and only 
if the dilated polytope $q{\Cal P}(D)$ is non-empty. By Theorem 6.3.5 of \cite{BR91} there 
exists an $n \times n$ symmetric non-negative integer matrix with zero trace, row/column sums 
$qd_1, \ldots, qd_n$ and the entries not exceeding $q$ if and only if 
$$\sum_{i=1}^k qd_i \ \leq \ qk(k-1) + \sum_{i=k+1}^n \min \left\{ qk, qd_i\right\} 
\quad k=1, \ldots, n. $$
Hence if the Erd\H{o}s-Gallai conditions are satisfied, the polytope $q{\Cal P}(D)$ is non-empty, and 
hence the polytope ${\Cal P}(D)$ is non-empty. 
{\hfill \hfill \hfill} \qed
\enddemo 

Next, we prove a sufficient condition for the polytope ${\Cal P}(D)$ to have a non-empty interior.

\proclaim{(12.2) Lemma}  Let $D=\left(d_1, \ldots, d_n \right)$ be a sequence of positive integers
such that 
$$d_1 \geq \ldots \geq d_n$$
and the strict Erd\H{o}s-Gallai conditions 
$$\sum_{i=1}^k d_i \ < \ k(k-1) + \sum_{i=k+1}^n \min\left\{k, d_i \right\} 
\quad \text{for} \quad k=1, \ldots, n$$
are satisfied. Then ${\Cal P}(D)$ has a non-empty interior.
\endproclaim
\demo{Proof} For a sufficiently small rational $\epsilon \geq 0$, let us define 
$$d_i(\epsilon)={d_i-(n-1) \epsilon \over 1-2\epsilon} \quad \text{for} \quad i=1, \ldots, n.$$
Clearly, $d_i(0)=d_i$ and 
$$d_1(\epsilon) \geq \ldots \geq d_n(\epsilon).$$
For all sufficiently small $\epsilon >0$ we have $d_i(\epsilon) >0$ for 
$i=1, \ldots, n$ and the Erd\H{o}s-Gallai conditions of Lemma 12.1 are satisfied for $d_i(\epsilon)$. 
Therefore, by Lemma 12.1, the polytope ${\Cal P}_{\epsilon}={\Cal P}\left(D_{\epsilon}\right)$ 
for $D=\bigl(d_1(\epsilon), \ldots, d_n(\epsilon)\bigr)$ is non-empty. Let $x=\left(\xi_{\{j, k\}}\right)$,
$x \in {\Cal P}_{\epsilon}$, be a point. Then the point $y=\left(\eta_{\{j, k\}} \right)$ 
defined by 
$$\eta_{\{j, k\}} = (1-2\epsilon) \xi_{\{j, k\}} + \epsilon \quad \text{for all} \quad 1 \leq j \ne k \leq n$$
is the desired interior point in ${\Cal P}(D)$.
{\hfill \hfill \hfill} \qed
\enddemo 
 
Next, we prove that our conditions on the minimum and maximum degree ensure that 
${\Cal P}(D)$ has a non-empty interior.

\proclaim{(12.3) Lemma} Let us fix real numbers $0 < \alpha < \beta < 1$ such that 
$$\beta \ < \ 2 \sqrt{\alpha} - \alpha, \quad \text{or, equivalently,} 
\quad (\alpha+\beta)^2 < 4 \alpha.$$ and let 
$D=\left(d_1, \ldots, d_n \right)$ be an integer sequence such that 
$$\alpha \ < \ {d_i \over n-1} \ < \ \beta \quad \text{for} \quad i=1, \ldots, n.$$
Then for 
$$n \ > \ \max \left\{ {\beta \over \alpha(1-\beta)}, \ {4(\beta-\alpha) \over 4 \alpha - (\alpha + \beta)^2} \right\}+1$$
the polytope ${\Cal P}(D)$ has a non-empty interior.
\endproclaim 
\demo{Proof} Without loss of generality, we assume that 
$$d_1 \geq \ldots \geq d_n.$$
Let us show that the strict Erd\H{o}s-Gallai conditions 
$$\sum_{i=1}^k d_i \ < \ k(k-1) + \sum_{i=k+1}^n \min\left\{k, d_i \right\} \quad 
\text{for} \quad k=1, \ldots, n$$
are satisfied.

We consider three different cases for $k$. 

Suppose that $k \leq \alpha (n-1)$. 
Then 
$$\sum_{i=1}^k d_i \ < \ k \beta(n-1) \quad 
\text{and} \quad \min\left\{k, d_i \right\}=k \quad \text{for all} \quad i.$$ 
Therefore,
$$k(k-1) + \sum_{i=k+1}^n \min\left\{k, d_i\right\} = k(k-1) + k(n-k)=k(n-1)$$
and the strict Erd\H{o}s-Gallai conditions are satisfied.

Suppose that $k \geq \beta(n-1)$. 
Then 
$$\sum_{i=1}^k d_i \ < \ k \beta(n-1) \quad \text{and} \quad \min\left\{k, d_i \right\} =d_i \ > \ 
\alpha(n-1).$$
If $k \geq \beta(n-1)+1$ then 
$k(k-1) \ \geq \ k \beta(n-1)$
and the strict Erd\H{o}s-Gallai conditions are satisfied. If $\beta(n-1) \leq k \leq \beta(n-1)+1$ then 
$$\split k(k-1) + \sum_{i=k+1}^n \min\left\{k, d_i \right\} \ \geq \ &\beta(n-1)(k-1) + \alpha(n-1)(n-k)\\
=\ &\beta k(n-1) +(n-1)\Bigl(\alpha(n-k)-\beta\Bigr) \\ 
\geq \ & \beta k(n-1) + (n-1) \Bigl(\alpha (n-1) (1-\beta)  - \beta \Bigr)\endsplit $$
and the strict Erd\H{o}s-Gallai conditions are satisfied provided 
$$n \ \geq \ { \beta \over \alpha(1-\beta)}+1.$$

Finally, suppose that $k=\gamma(n-1)$ for some $\alpha < \gamma < \beta$. 
Then 
$$\sum_{i=1}^k d_i \ < \ k \beta(n-1)=\gamma\beta(n-1)^2 \quad \text{and} \quad \min\left\{k, d_i \right\} > \ \alpha(n-1).$$
Therefore,
$$k(k-1) + \sum_{i=k+1}^n \min\left\{k, d_i \right\} \ > \ 
\gamma^2(n-1)^2-(\gamma-\alpha)(n-1) + (1-\gamma)\alpha(n-1)^2.$$
The minimum value of the function 
$$\gamma \longmapsto \gamma^2 + (1-\gamma) \alpha - \gamma \beta$$
is attained at $\gamma=(\alpha +\beta)/2$ and equal to 
$$\alpha-{(\alpha+\beta)^2 \over 4} > 0$$
Therefore, the strict Erd\H{o}s-Gallai conditions are satisfied, provided 
$$n \ > \ {4(\beta-\alpha) \over 4 \alpha - (\alpha + \beta)^2} +1.$$
The proof now follows by Lemma 12.2.
{\hfill \hfill \hfill} \qed
\enddemo

\subhead (12.4) Proof of Theorem 2.1 \endsubhead  By Lemma 12.3, the polytope ${\Cal P}(D)$ 
contains a point $y=\left(\eta_{\{j, k\}}\right)$ such that 
$0 < \eta_{\{j, k\}} <1$ for all $j, k$.
 First, we show that the maximum entropy matrix $z$ lies in the 
 interior of ${\Cal P}(D)$, that is, $0 < \zeta_{\{j, k\}} < 1$ for all $j, k$. 
 
 We have 
 $${\partial \over \partial \xi_{\{j, k\}}} H(x)=\ln {1 -\xi_{\{j, k\}} \over \xi_{\{j, k\}}}.$$
 We note that the value of the derivative is $+\infty$ at $\xi_{\{j, k\}}=0$ (we consider the right derivative 
 there), is $-\infty$ at $\xi_{\{j, k\}}=1$ (we consider the left derivative there) and is finite for 
 $0 < \xi_{\{j, k\}} <1$. Therefore, if for the maximum point $z$ and some $j \ne k$ we have 
 $\zeta_{\{j, k\}} \in \{0, 1\}$ then for
 $\tilde{z}=(1-\epsilon) z  + \epsilon y$ for a sufficiently small $\epsilon >0$, we have 
 $\tilde{z} \in {\Cal P}(D)$ and $H\left(\tilde{z}\right) > H(z)$, which is a contradiction. 
 
 Since the maximum value of $H$ is attained at an interior point of ${\Cal P}(D)$, the gradient of 
 $H$ at the maximum point is orthogonal to the affine span of ${\Cal P}(D)$, that is, 
 $$\ln {1 -\zeta_{\{j, k\}} \over \zeta_{\{j, k\}}} = \lambda_j + \lambda_k,$$
 or, equivalently,
 $$\zeta_{\{j, k\}} = {1 \over 1 +e^{\lambda_j +\lambda_k}} \quad \text{for all} \quad 1 \leq j \ne k \leq n\tag12.4.1$$ 
 for some real $\lambda_1, \ldots, \lambda_n$. 
 Without loss of generality, we assume that 
  $$\lambda_1 \ \leq \ \lambda_j \ \leq \ \lambda_n \quad \text{for all} \quad j. \tag12.4.2$$
 From the choice of $\epsilon$ in Theorem 2.1, it follows that 
 $$\epsilon \ \leq \ \alpha \quad \text{and} \quad \beta \ \leq \ 2 \sqrt{\alpha -\epsilon} - \alpha.
  \tag12.4.3$$
 Our next goal is to show that 
 $$\lambda_n \ \leq \ 2\ln {1 \over \epsilon}. \tag12.4.4$$
 Aiming for a contradiction, suppose that 
 $$\lambda_n \ > \ 2\ln {1 \over \epsilon}. $$
 Then, necessarily, 
 $$\lambda_1\ < \ \ln \epsilon$$ 
since otherwise by (12.4.1) and (12.4.2) we have 
 $$\zeta_{\{j, n\}} = {1 \over 1 + e^{\lambda_j + \lambda_n}} \ \leq \ {1 \over 1+ e^{\lambda_1 +\lambda_n}} \ < \ \epsilon$$ 
 and 
 $$d_n=\sum_{j: \ j \ne n} \zeta_{\{j, n\}}  \ < \ \epsilon (n-1),$$
 which by (12.4.3) contradicts the lower bound for $d_i$.
 
Since $\lambda_1 < \ln \epsilon$ and $\lambda_n > -2\ln \epsilon$, we 
deduce from (12.4.1) that for $1 < j < n$ we have 
$$\aligned \zeta_{\{j, n\}} = &{1 \over 1+ e^{\lambda_j + \lambda_n}} \ < \ {1 \over 1 + e^{\lambda_n}} \ < \epsilon 
\quad \text{provided} \quad  \lambda_j \geq 0 \\ &\qquad \qquad \text{and} \\
\zeta_{\{1, j\}}= &{1 \over 1 + e^{\lambda_1 + \lambda_j}} \ > \ {1 \over 1+ e^{\lambda_1}} 
\ > 1-\epsilon \quad \text{provided} \quad \lambda_j \leq 0. \endaligned \tag12.4.5  $$
Denoting 
$$\tau=\zeta_{\{1, n\}} = {1 \over 1 + e^{\lambda_1 +\lambda_n}} < 1,$$
by (12.4.1) and (12.4.2) we obtain 
that for $1< j < n$ we have 
$$\aligned \zeta_{\{1, j\}} =&{1 \over 1+ e^{\lambda_1 + \lambda_j}} \ \geq \ {1 \over 1 + e^{\lambda_1 + \lambda_n}} =\tau \\ 
&\qquad \qquad \text{and} \\ 
\zeta_{\{j, n\}} =&{1 \over 1+e^{\lambda_j + \lambda_n}} \ \leq \ {1 \over 1+ e^{\lambda_1 +\lambda_n}}
=\tau.
\endaligned \tag12.4.6$$
Let 
$$\left|\bigl\{ 1\leq j < n: \quad \lambda_j \leq 0 \bigr\}\right| = \gamma (n-1) \quad \text{for some} \quad 
0 \leq \gamma \leq 1.$$
Combining (12.4.5) and (12.4.6), we obtain 
$$\split \beta (n-1) \ > \ &d_1 =\sum_{j: \ j \ne 1} \zeta_{\{1, j\}} = 
\sum_{j\ne 1:  \ \lambda_j  \leq 0} \zeta_{\{1, j\}} + \sum_{j \ne 1: \ \lambda_j > 0} \zeta_{\{1, j\}} 
\\  > \ &(1-\epsilon) \gamma (n-1) + (n-1)(1-\gamma) \tau \endsplit$$
and 
$$\split \alpha (n-1) \ < \ &d_n =\sum_{j: \ j \ne n} \zeta_{\{j, n\}} = 
\sum_{j\ne n:  \ \lambda_j >0} \zeta_{\{j, n\}} + \sum_{j \ne n: \ \lambda_j \leq  0} \zeta_{\{j, n\}} 
\\  \leq \ &\epsilon (1- \gamma) (n-1) + (n-1)\gamma \tau. \endsplit$$
 Consequently,
 $$\beta \ > \ (1-\epsilon)\gamma + (1-\gamma) \tau \quad \text{and} \quad \alpha 
 \ < \ \epsilon(1-\gamma) + \gamma \tau.$$
 Therefore,
 $$\beta+ \epsilon \ > \ \gamma + (1-\gamma) \tau \quad \text{and} \quad 
 \alpha -\epsilon \ < \ \gamma \tau.$$
 Since the function $2\sqrt{x} -x$ is increasing for $0 < x <1$, from (12.4.3)
 it follows that 
 $$\gamma + (1-\gamma) \tau \ < \ 2\sqrt{\gamma \tau} - \gamma \tau,$$
 or, equivalently, 
 $${\gamma + \tau \over 2} \ < \ \sqrt{\gamma \tau},$$
 which is a contradiction.
 
 The contradiction shows that (12.4.4) indeed holds. Then, by (12.4.2), we have 
 $$\lambda_j \ \leq \ 2 \ln {1 \over \epsilon} \quad \text{for} \quad j=1, \ldots, n.$$
We claim now that 
$$\lambda_1 \ \geq \ 3 \ln \epsilon. \tag12.4.7$$
Indeed, if $\lambda_1 < 3 \ln \epsilon$ 
then by (12.4.1)
$$\zeta_{\{1, j\}} = {1 \over 1+ e^{\lambda_1 + \lambda_j}} \ > \ {1 \over 1+ \epsilon} \ > \ 1-\epsilon
\quad \text{for} \quad j=1, \ldots, n-1$$
and 
$$\beta(n-1) \ > \ d_1 = \sum_{j: \ j \ne n} \zeta_{\{1, j\}} \ > (n-1)(1-\epsilon),$$
which contradicts (12.4.3).

Summarizing, from (12.4.4) and (12.4.7), we obtain 
$${\epsilon^4 \over 1+\epsilon^4} \ \leq \ {1 \over 1 + e^{2 \lambda_n}} \ \leq \ \zeta_{\{j, k\}} \ \leq \ {1 \over 1+ e^{2 \lambda_1}} \ \leq\ 
{1 \over 1 +\epsilon^6} \quad \text{for all} \quad j \ne k,$$
which completes the proof.
{\hfill \hfill \hfill} \qed

\head Acknowledgment \endhead

The authors are grateful to Brendan McKay and Sourav Chatterjee for useful comments.

\enddocument
\end